\documentclass[a4paper, 12pt]{article}
\usepackage{amsmath,amsfonts,amssymb}
\usepackage[all]{xy}
\def\comack{\mathsf{M}^c}

\def\per{\mathsf{perm}}
\def\fun{\mathsf{Fun}}
\def\iiota{\imath}
\def\jiota{\jmath}
\newcommand{\bn}[2]{\displaystyle\binom{#1}{#2}}

\makeatletter
\newbox\bk@bxb
\newbox\bk@bxa
\newif\if@bkcont
\newcount\bk@lcnt

\def\breakboxskip{2pt}
\def\breakboxparindent{1.8em}

\def\breakbox{\vskip\breakboxskip\relax
\setbox\bk@bxb\vbox\bgroup
\advance\linewidth -2\fboxrule
\hsize\linewidth\@parboxrestore
\parindent\breakboxparindent\relax}

\def\bk@split{%
\@tempdimb\ht\bk@bxb % height of original box
\advance\@tempdimb\dp\bk@bxb
\setbox\bk@bxa\vsplit\bk@bxb to\z@ % split it
\setbox\bk@bxa\vbox{\unvbox\bk@bxa}% recover height & depth of \bk@bxa
\setbox\@tempboxa\vbox{\copy\bk@bxa\copy\bk@bxb}% naive concatenation
\advance\@tempdimb-\ht\@tempboxa
\advance\@tempdimb-\dp\@tempboxa}% gap between two boxes

\def\bk@addfsepht{%
\setbox\bk@bxa\vbox{\vskip\fboxsep\box\bk@bxa}}

\def\bk@addskipht{%
\setbox\bk@bxa\vbox{\vskip\@tempdimb\box\bk@bxa}}

\def\bk@addfsepdp{%
\@tempdima\dp\bk@bxa
\advance\@tempdima\fboxsep
\dp\bk@bxa\@tempdima}

\def\bk@addskipdp{%
\@tempdima\dp\bk@bxa
\advance\@tempdima\@tempdimb
\dp\bk@bxa\@tempdima}

\def\bk@line{%
\hbox to \linewidth{%
\hskip-2\fboxsep\vrule \@width\fboxrule\hskip.5\fboxsep\vrule \@width\fboxrule\hskip1.5\fboxsep
\box\bk@bxa\hfil
}}%

\def\endbreakbox{\egroup
\ifhmode\par\fi{\noindent\bk@lcnt\@ne
\@bkconttrue\baselineskip\z@\lineskiplimit\z@
\lineskip\z@\vfuzz\maxdimen
\bk@split\bk@addfsepht\bk@addskipdp
\ifvoid\bk@bxb % Only one line
\def\bk@fstln{\bk@addfsepdp
\hskip-\parindent\vbox{\llap{\raisebox{-2ex}{\rule{1.5\fboxsep}{\fboxrule}\hskip.5\fboxsep}}\bk@line\llap{\rule{1.5\fboxsep}{\fboxrule}\hskip.5\fboxsep}}}% ajouté le saut négatif...

\else % More than one line
\def\bk@fstln{\vbox{\llap{\raisebox{-2ex}{\rule{1.5\fboxsep}{\fboxrule}\hskip.5\fboxsep}}\bk@line}\hfil%
\advance\bk@lcnt\@ne
\loop
\bk@split\bk@addskipdp\leavevmode
\ifvoid\bk@bxb % The last line
\@bkcontfalse\bk@addfsepdp
\vtop{\bk@line\noindent\hskip-2\fboxsep{\rule{1.5\fboxsep}{\fboxrule}}}%
\else % 2,...,(n-1)
\bk@line
\fi
\hfil\advance\bk@lcnt\@ne
\if@bkcont\repeat}%
\fi
\leavevmode\bk@fstln\par}\vskip\breakboxskip\relax}

\def\smp{\smallskip\par}
\def\un{{\bf 1}}
\def\zero{\{0\}}
\def\pf{\noindent{\bf Proof~:}\ }
\def\findemo{~\leaders\hbox to 1em{\hss\  \hss}\hfill~\raisebox{.5ex}{\framebox[1ex]{}}\smp}
\def\spn{\bigskip\par\noindent}
\def\mpn{\medskip\par\noindent}
\def\smpn{\smallskip\par\noindent}
\def\normal{\mathop{\trianglelefteq}}
\def\sp{\bigskip\par}
\def\smp{\smallskip\par}
\def\smpn{\smallskip\par\noindent}

\def\mpoint{\;\;.}
\def\mvirg{\;\;,}

\def\Res{{\rm Res}}
\def\Ind{{\rm Ind}}
\def\Inf{{\rm Inf}}

\def\Iso{{\rm Iso}}

\def\Hom{{\rm Hom}}
\def\End{{\rm End}}
\def\Ext{{\rm Ext}}
\def\Inf{{\rm Inf}}

\def\Ker{{\rm Ker}}
\def\Coker{{\rm Coker}}
\def\Id{{\rm Id}}

\def\dsp{\displaystyle}
\def\Z{\mathbb{Z}}
\def\N{\mathbb{N}}
\def\F{\mathbb{F}}

\newcommand{\dirsum}[1]{\mathop{\oplus}_{#1}\limits}
\newcommand{\romain}[1]{\uppercase\expandafter{\romannumeral #1}}

\newcommand{\flh}[2]{\mathop{\hbox to 12mm{\rightarrowfill}}_{\displaystyle #2}^{\displaystyle #1}\limits}
\newcommand{\sflh}[2]{\mathop{\hbox to 12mm{\rightarrowfill}}_{\scriptstyle #2}^{\scriptstyle #1}\limits}

\newcommand{\gMod}[1]{#1{\hbox{-}\mathsf{Mod}}}
\newcommand{\gmod}[1]{#1\hbox{-$\mathsf{mod}$}}

\newcommand{\sur}[1]{\,\overline{\! #1}}

\newcommand{\oplusb}[2]{\mathop{\bigoplus}_{{\scriptstyle #1}\atop{\scriptstyle #2}}}

\newcommand{\carre}[8]{\begin{array}{ccc}
#1&\mathop{\hbox to 12mm{\rightarrowfill}}^{\displaystyle{#2}}\limits&#3\\
\llap{$\displaystyle{#4}$}\left\downarrow\vbox to 6mm{}\right. & & \left\downarrow\vbox to 6mm{}\right.\rlap{$\displaystyle{#5}$}\\
#6&\mathop{\hbox to 12mm{\rightarrowfill}}_{\displaystyle #7}\limits&#8\\
\end{array}}

\newcommand{\carrem}[8]{\begin{array}{ccc}
#1&\mathop{\hbox to 12mm{\rightarrowfill}}^{\displaystyle #2}\limits&#3\\
\llap{$\displaystyle #4$}\left\uparrow\vbox to 6mm{}\right. & & \left\uparrow\vbox to 6mm{}\right.\rlap{$\displaystyle #5$}\\
#6&\mathop{\hbox to 12mm{\rightarrowfill}}_{\displaystyle #7}\limits&#8\\
\end{array}}

\newenvironment{enonce}[1]{\pagebreak[2]\refstepcounter{subsection}\refstepcounter{prop}\smpn{{\bf \thesection.\arabic{prop}.\ \ #1~:}}\begin{it} }{\end{it}\smp}
\newenvironment{enonce*}[1]{\pagebreak[2]\smpn{#1~:}\begin{it} }{\end{it}\smp}
\newcommand{\result}[1]{\begin{enonce}{#1}}
\def\fresult{\end{enonce}}
\newcommand{\npar}{\smallskip\par\noindent\pagebreak[2]\refstepcounter{subsection}\refstepcounter{prop}{\bf \thesection.\arabic{prop}.\ \ }}

\newenvironment{mth}[1]{\begin{breakbox}\begin{enonce}{#1}}{\end{enonce}\end{breakbox}}
\newenvironment{mth*}[1]{\begin{breakbox}\begin{enonce*}{#1}}{\end{enonce*}\end{breakbox}}
\newenvironment{rem}[1]{\refstepcounter{subsection}\refstepcounter{prop} \mpn{{\bf \thesection.\arabic{prop}.}\ \ \bf#1\ :}}{\smp}

\def\dom{\backslash}
\makeatletter
\renewenvironment{enumerate}{\ifnum \@enumdepth >3 \@toodeep\else
      \advance\@enumdepth \@ne
      \edef\@enumctr{enum\romannumeral\the\@enumdepth}\list
      {\csname label\@enumctr\endcsname}{\setlength{\topsep}{1ex}\setlength{\itemsep}{0pt}\usecounter
        {\@enumctr}\def\makelabel##1{\hss\llap{##1}}}\fi}{\endlist}
\renewenvironment{itemize}{\ifnum \@itemdepth >3 \@toodeep\else \advance\@itemdepth \@ne
\edef\@itemitem{labelitem\romannumeral\the\@itemdepth}%
\list{\csname\@itemitem\endcsname}{\setlength{\topsep}{1ex}\setlength{\itemsep}{0pt}\def\makelabel##1{\hss\llap{##1}}}\fi}
{\endlist}
\makeatother
\makeatletter
\def\@sect#1#2#3#4#5#6[#7]#8{\ifnum #2>\c@secnumdepth
    \let\@svsec\@empty\else
    \refstepcounter{#1}\edef\@svsec{\csname the#1\endcsname .\hskip .5em}\fi
    \@tempskipa #5\relax
     \ifdim \@tempskipa>\z@
       \begingroup #6\relax
         \@hangfrom{\hskip #3\relax\@svsec}{\interlinepenalty \@M #8\par}%
       \endgroup
      \csname #1mark\endcsname{#7}\addcontentsline
        {toc}{#1}{\ifnum #2>\c@secnumdepth \else
                     \protect\numberline{\csname the#1\endcsname}\fi
                   #7}\else
       \def\@svsechd{#6\hskip #3\relax  %% \relax added 2 May 90
                  \@svsec #8\csname #1mark\endcsname
                     {#7}\addcontentsline
                          {toc}{#1}{\ifnum #2>\c@secnumdepth \else
                            \protect\numberline{\csname the#1\endcsname}\fi
                      #7}}\fi
    \@xsect{#5}}
\def\section{\@startsection {section}{1}{\z@}{-3.5ex plus-1ex minus
    -.2ex}{2.3ex plus.2ex}{\reset@font\Large\bf}}  %{\reset@font\large\bf}}

\makeatother
\renewenvironment{equation}{\refstepcounter{subsection}\refstepcounter{prop}$$}{\leqno{\bf (\theprop)}$$}

\def\mar[#1]{\ar@{-}[#1]|-{\object@{<}}}
\def\marb[#1]{\ar@{-}[#1]|{\object+{  }}}

\def\Ker{{\rm Ker}}
\def\Coker{{\rm Coker}}
\def\Img{{\rm Im}}

\def\Hom{{\rm Hom}}
\def\Ext{{\rm Ext}}
\def\Res{{\rm Res}}

\def\tr{{\rm tr}}
\newcommand{\bin}[3]{\displaystyle{\binom{S_#1^#3}{S_#2^#3}}}
\def\F{{\mathbb{F}}}
\def\A{{\mathcal A}}

\def\Aa{{\mathcal A}}

\def\Ee{{\mathcal E}}
\def\Oo{{\mathcal O}}

\def\rk{{\rm pos}\,}

\begin{document}

\title{\bf The extension algebra\\of some cohomological Mackey functors}
\author{Serge Bouc and Radu Stancu}
\date{\ }
\maketitle
\noindent{\footnotesize {\bf Abstract~:} Let $k$ be a field of characteristic $p$. We construct a new inflation functor for cohomological Mackey functors for finite groups over $k$. Using this inflation functor, we give an explicit presentation of the graded algebra of self extensions of the simple functor $S_\un^G$, when $p$ is odd and $G$ is an elementary abelian $p$-group. \spn
{\bf Keywords~:} Cohomological, Mackey functor, extension, simple.}
%%%%%%%%%%%%%%%%%%%%%%%%%
\section{Introduction}%%%
%%%%%%%%%%%%%%%%%%%%%%%%%
Let $k$ be a field and $G$ be a finite group. The theory of {\em Mackey functors} and {\em cohomological Mackey functors} for $G$ over $k$ originates in the work of Green~(\cite{jagreen}) and Dress~(\cite{dress}), at the beginning of the 70's. It can be viewed as the theory of induction and restriction, when we forget the particular framework of linear representations of $G$ over $k$. Many important developments have been achieved since, culminating in the comprehensive and seminal paper by Th\'evenaz-Webb (\cite{thevwebb}) in 1995, where the authors introduce {\em the Mackey algebra} $\mu_k(G)$, and show, among many other fundamental results, that the category of Mackey functors for $G$ over $k$ is equivalent to the category of $\mu_k(G)$-modules. Similarly, they show that the subcategory $\comack_k(G)$ of cohomological Mackey functors for $G$ over $k$ is equivalent to the category of $co\mu_k(G)$-modules, where $co\mu_k(G)$ is a specific quotient of $\mu_k(G)$, called {\em 
 the cohomological Mackey algebra}.\par
The algebras $\mu_k(G)$ and $co\mu_k(G)$ share many similarities with the group algebra $kG$~: e.g., they are finite dimensional $k$-vector spaces, of dimension independent on $k$, the Maschke theorem holds, there is a good theory of decomposition from characteristic 0 to characteristic $p$, etc\ldots~These resemblances raise some natural questions, whether a given theorem on $kG$ will admit an analogue for $\mu_k(G)$ or $co\mu_k(G)$.\par
This was the main motivation in~\cite{cohocplx}, where the question of complexity of cohomological Mackey functors was solved (in the only non-trivial case where $k$ is a field of positive characteristic $p$ dividing the order of $G$). It was also shown there how this question can be reduced to the consideration of elementary abelian $p$-groups $E$ appearing as subquotients of $G$, and to the knowledge of enough information on the algebra $\Ee=\Ext^*(S_\un^E,S_\un^E)$ of self extensions of a particular simple functor $S_\un^E$ for these groups. Along the way, a presentation of this algebra was given when $p=2$, together with a formula for the Poincar\'e series. In the case $p>2$, no such presentation was given, and a conjecture was proposed for the Poincar\'e series of $\Ee$. This conjecture was only proved in the case $p=3$.\par
This paper settles completely the case $p>2$~: a presentation of $\Ee$ is given, and, as a corollary, the forementioned conjecture is proved. The main results are the following, where $S_{1,W}^H$ denotes the simple functor for the group $H$ defined as in~\ref{definition simple functor}. To simplify notation, when $W=k$ is the trivial module, we drop this subscript. We start by the construction of a new inflation functor for cohomological Mackey functors~: 
\begin{mth}{Theorem} \label{inflation functor} Let $k$ be a field, let $G$ be a finite group, let $N\normal G$ and let $V$ be a simple $k(G/N)$-module. Then there exist an exact functor $\sigma_{G/N}^G$ from $\comack_k(G/N)$ to $\comack_k(G)$ satisfying
$$\sigma_{G/N}^G(S_{\un,V}^{G/N})=S_{\un,\Inf_{G/N}^GV}^G\mpoint$$
Suppose, moreover, that $G=N\rtimes H$ is the semidirect product of $N$ by a group~$H$.
\begin{enumerate}
\item If $V$ is a $kH$-module, let $\tilde{V}$ be the $kG$-module $\Inf_{G/N}^G\Iso_H^{G/N}V$. Then the restriction of $\tilde{V}$ to $H$ is isomorphic to $V$.
\item The composition $\Res_H^G\sigma_{G/N}^G\Iso_{H}^{G/N}$ is isomorphic to the identity functor of $\comack_k(H)$.
\item Let $V$ and $W$ be simple $kH$-modules. Then, for any $n\in \N$, the restriction from $G$ to $H$ induces a split surjection
$$r_H^G:\Ext^n_{\comack_k(G)}(S_{\un,\tilde{V}}^G,S_{\un,\tilde{W}}^G)\to\Ext^n_{\comack_k(H)}(S_{\un,V}^H,S_{\un,W}^H)\mpoint$$
\end{enumerate}
\end{mth}
Let now $p$ be an odd prime, $k$ be a field of characteristic $p$ and $G$ be an elementary abelian $p$-group of rank $r$. Then we can give an explicit presentation of $\Ext^*_{\comack_k(G)}(S_\un^G,S_\un^G)$, the graded algebra of self extensions of the simple functor~$S_\un^G=S_{\un,k}^G$. Here is the main result of this paper~:
\begin{mth}{Theorem}\label{presentation Ext algebra 1}
Let $\A$ be the graded $k$-algebra with generators 
$$\{\hat\tau_\varphi\mid\varphi\in\Hom(G,\F_p^+)\},\;\hbox{in degree 1, and }\{\hat\gamma_X|X\le G,\, |X|=p\},\;\hbox{in degree 2},$$ 
subject to the relations 
\begin{itemize}
\item [(L1)] $\forall \varphi, \psi\in \Hom(G,\F_p^+),\;\;\hat\tau_{\varphi +\psi}=\hat\tau_\varphi+\hat\tau_\psi$.
\item [(L2)] If $p\geq 5$, then $\forall \varphi\in \Hom(G,\F_p^+),\;\;\hat\tau_\varphi^2=0$, and $[\hat\tau_\varphi,\sum_{X\nleq \Ker\varphi}\limits\hat\gamma_X]=0$. \\
If $p=3$, then $\hat\tau_\varphi^2=-\sum_{X\nleq \Ker\varphi}\limits\hat\gamma_X$.
\item [(L3)] $\forall \varphi\in\Hom(G,\F_p^+), \forall X,\;|X|=p,\;X\leq \Ker\varphi,\;\;[\hat\tau_\varphi,\gamma_X]=0$.
\item [(L4)] $\forall Q\leq G, |Q|=p^2, \forall X<Q, |X|=p,\;\;[\hat\gamma_X,\sum_{Y<Q}\limits\hat\gamma_Y]=0$.
\end{itemize}
Then there is an isomorphism of graded algebras from $\A$ to $\Ext^*_{\comack_k(G)}(S_\un^G,S_\un^G)$.
\end{mth}
In view of $(L1)$, this set of generators is redundant~: we choose a direct sum decomposition $G=Y_1\oplus Y_2\oplus\cdots\oplus Y_{r-1}\oplus Y_r$, and for $1\leq i\leq r$, we choose a group homomorphism $\varphi_i:G\to\F_p$ with kernel $\oplusb{1\le j\le r}{j\ne i}\limits Y_j$. 
\begin{mth}{Theorem}\label{presentation Ext algebra}
Let $\A$ be the graded $k$-algebra with generators 
$$\{\hat\tau_i|1\le i\le r\}\;\hbox{in degree 1, and }\{\hat\gamma_X|X\le G,\, |X|=p\}\;\hbox{in degree 2},$$ 
subject to the relations 
\begin{itemize}
\item [(R1)]$\hat\tau_i\hat\tau_i=0$, if $p\geq 5$, or \\
$\hat\tau_i\hat\tau_i=-\sum_{X\nleq\Ker\,\varphi_i}\limits\hat\gamma_X$ if $p=3$.
\item [(R2)]$\hat\tau_i\hat\tau_j+\hat\tau_j\hat\tau_i=0$ for $1\le j\le i \le r$, if $p\geq 5$, or \\
$\hat\tau_i\hat\tau_j+\hat\tau_j\hat\tau_i=\sum_{X\nleq\Ker(\varphi_i+\varphi_j)}\limits\hat\gamma_X-\sum_{X\nleq\Ker\,\varphi_i}\limits\hat\gamma_X-\sum_{X\nleq\Ker\,\varphi_j}\limits\hat\gamma_X$ if $p=3$.
\item [(R3)]$[\hat\tau_i,\sum_{X\nleq\Ker\,\varphi_i}\limits\hat\gamma_X]=0$;
\item [(R4)]$[\varphi_j(x)\hat\tau_i-\varphi_i(x)\hat\tau_j,\hat\gamma_{\langle x\rangle}]=0, \text{ for }1\le i<j\le 1,\,x\in G$;
\item [(R5)]$[\hat\gamma_X,\sum_{Y<Q}\limits\hat\gamma_Y]=0,\text{ for all }\,X<Q\le G,\,|X|=p,\,|Q|=p^2$.
\end{itemize}
Then there is an isomorphism of graded algebras from $\A$ to $\Ext^*_{\comack_k(G)}(S_\un^G,S_\un^G)$.
\end{mth}
As a consequence of the above theorem, we get the following result, which was proved for $p=3$ in (\cite{cohocplx},~Theorem~14.2)~ and conjectured for $p\ge 5$:
\begin{mth}{Proposition}\label{proven conjecture}Let $k$ be a field of odd characteristic $p$, and $G\cong (C_p)^r$. Then~:
\begin{enumerate}
\item The algebra $\Ee=\Ext^*_{\comack_k(G)}(S_\un^G,S_\un^G)$ is generated by the elements $\tau_\varphi^G$ in degree 1, for $\varphi\in \Hom(G,k^+)$, and by the elements $\gamma_X^G$ in degree~2, for $X\leq G$ with $|X|=p$.
\item The Poincar\'e series for $\Ee$ is equal to
$$\frac{1}{(1-t)\big(1-t-(p\!-\!1)t^2\big)\big(1-t-(p^2\!-\!1)t^2\big)\ldots\big(1-t-(p^{r-1}\!-\!1)t^2\big)}\mpoint$$
\end{enumerate}
\end{mth}
The paper is organized as follows. In Section 2 we construct a new inflation functor for cohomological Mackey functors for a finite group $G$ over a field $k$. In the case where the group is a semi-direct product, more properties of this functor are given. For completeness, Section 3 recalls a series of results on extensions of cohomological Mackey functors when $G$ is an elementary abelian $p$-group $G$ and $k$ is a field of characteristic $p$. These results, taken from \cite{cohocplx}, are crucial in the second part of the paper. The use of the inflation functor constructed in Section~2 makes the proof of the exactness of the sequence in Corollary~\ref{short exact sequence} straightforward. Section 4 deals with finding relations in the graded algebra $\Ee=\Ext^*_{\comack_k(G)}(S_\un^G,S_\un^G)$. Section 5 gives a recursive direct decomposition of this algebra. The decomposition is obtained through an involved induction on the rank of $G$, constructed on the technical support 
 of  Lemma~\ref{diagrchase}. Based on this direct decomposition and the relations stated in Section~4, Section~6 builds a presentation of $\Ee$. Some arithmetics of extensions in abelian categories, needed in the paper, are presented in the Appendix.

%%%%%%%%%%%%%%%%%%%%%%%%%%%%%%%%%%%%%%%%%%
\section{Yet another inflation functor}%%%
%%%%%%%%%%%%%%%%%%%%%%%%%%%%%%%%%%%%%%%%%%
\npar Let $k$ be a commutative ring with identity element, and $G$ be a finite group.  Let $\comack_k(G)$ denote the category of cohomological Mackey functors for $G$ over $k$.
By Yoshida's Theorem, the category $\comack_k(G)$ is equivalent to the category $\fun_k(G)$ of $k$-linear contravariant functors from the category $\per_k(G)$ of finitely generated permutation $kG$-modules, to the category $\gMod{k}$ of $k$-modules.
\npar When $G$ and $H$ are finite groups, any $k$-linear functor from $\per_k(G)$ to $\per_k(H)$ gives by precomposition a functor from $\fun_k(H)$ to $\fun_k(G)$. Two special cases of this situation have been considered in~\cite{cohocplx} (Section~3.12)~: if $U$ is a finite $(H,G)$-biset, the functor 
$$\mathsf{t}_U:V\mapsto kU\otimes_{kG}V$$
from $\per_k(G)$ to $\per_k(H)$ yields a functor
$$\mathsf{L}_U:F\mapsto F\circ \mathsf{t}_U$$
from $\fun_k(H)$ to $\fun_k(G)$. Similarly, the functor
$$\mathsf{h}_U:W\mapsto \Hom_{kH}(kU,W)$$
from $\per_k(H)$ to $\per_k(G)$ yields a functor 
$$\mathsf{R}_U:F\mapsto F\circ h_U$$
from $\fun_k(G)$ to $\fun_k(H)$, which is right adjoint to $\mathsf{L}_U$.
\npar Several special cases of this special case occur when $H=G/N$, where $N$ is a normal subgroup of $G$. First, if $U=G/N$, viewed as a $(G,H)$-biset in the obvious way, then the functor $\mathsf{L}_U:\fun_k(G/N)\to\fun_k(G)$ is denoted by $\iiota_{G/N}^G$, and $\mathsf{R}_U:\fun_k(G)\to\fun_k(G/N)$ is denoted by  $\rho_{G/N}^G$. \par
Reversing the actions, and viewing $V=G/N$ as an $(H,G)$-biset, the functor $\mathsf{R}_V:\fun_k(G/N)\to\fun_k(G)$ is denoted by $\jiota_{G/N}^G$, and $\mathsf{L}_V:\fun_k(G)\to\fun_k(G/N)$ again by $\rho_{G/N}^G$ (there is no notational conflict here, as $\mathsf{L}_V\cong\mathsf{R}_U$ by~\cite{cohocplx} Proposition~3.15).\par
The functors $\iiota_{G/N}^G$ and $\jiota_{G/N}^G$ are two inflation-like functors for cohomological Mackey functors. They are not isomorphic in general, and also both different from the usual inflation\footnote{Note that this inflation functor $\Inf_{G/N}^G$ does not preserve cohomological Mackey functors over $k$ unless $N$ is a $p$-group} $\Inf_{G/N}^G$ for Mackey functors (as defined in~\cite{thevwebb} Section~2). Similarly, the functor $\rho_{G/N}^G$ is a deflation functor.
\npar When $\mathcal{C}$ is a $k$-linear category, recall (\cite{sga4} Exemple 8.7.8 p. 97) that the {\em karoubian envelope} $\mathcal{C}^+$ of $\mathcal{C}$ is the category defined as follows~: the object of $\mathcal{C}^+$ are the pairs $(X,e)$, where $X$ is an object of $\mathcal{C}$ and $e$ is an idempotent in $\End_{\mathcal{C}}(X)$. If $(X,e)$ and $(Y,f)$ are objects of $\mathcal{C}^+$, then by definition
$$\Hom_{\mathcal{C}^+}\big((X,e),(Y,f)\big)=f\Hom_{\mathcal{C}}(X,Y)e\mvirg$$
and the composition of morphism in $\mathcal{C}^+$ is induced by the composition of morphism in $\mathcal{C}$.\par
The category $\mathcal{C}^+$ is a $k$-linear category. Moreover, the correspondence 
$$\left\{\begin{array}{ccl}X\in\mathcal{C}&\mapsto& (X,\Id_X)\in\mathcal{C}^+\\f\in\Hom_{\mathcal{C}}(X,Y)&\mapsto& f\in\Hom_{\mathcal{C}^+}\big((X,\Id_X),(Y,\Id_Y)\big)\end{array}\right.$$
is a fully faithful $k$-linear functor $\mathsf{i}$ from $\mathcal{C}$ to $\mathcal{C}^+$.\par
By composition, this functor induces a functor $\mathsf{I}:M\mapsto M\circ\mathsf{i}$ from the category $\mathcal{F}^+_k$ of $k$-linear functors from $\mathcal{C}^+$ to $\gMod{k}$, to the category $\mathcal{F}_k$ of $k$-linear functors from $\mathcal{C}$ to $\gMod{k}$, which is easily seen to be an equivalence of categories~: the functor $\mathsf{J}:\mathcal{F}_k\to \mathcal{F}^+_k$ defined by
$$\mathsf{J}(L)\big(X,e)\big)=L(e)\big(L(X)\big)$$
is a quasi-inverse to $\mathsf{I}$.
\npar Applying this to the category $\mathcal{C}=\per_k(G)$, we observe that the category $\mathcal{C}^+$ is equivalent to the full subcategory $\per^+_k(G)$ of $\gmod{kG}$ consisting of direct summands of finitely generated permutation $kG$-modules~: this equivalence is induced by the functor $(X,e)\mapsto e(X)$. It follows that the category $\comack_k(G)$ is equivalent to the category $\fun^+_k(G)$ of contravariant $k$-linear functors from $\per^+_k(G)$ to $\gMod{k}$. 
\npar It follows that when $G$ and $H$ are finite groups, any $k$-linear functor from $\per^+_k(G)$ to $\per^+_k(H)$ induces by composition a functor from $\fun^+_k(H)$ to $\fun^+_k(G)$, and by the above remarks, this yields a corresponding functor $\comack_k(H)\to \comack_k(G)$. In particular, when $X$ is a direct summand of a finitely generated permutation $(kH,kG)$-bimodule, the functor $t_X=X\otimes_{kG}{-}$ yields a functor $L_X:\comack_k(H)\to \comack_k(G)$, and the functor $h_X=\Hom_{kH}(X,{-})$ yields a functor $R_X:\comack_k(G)\to \comack_k(H)$.
\npar We will consider here another special case of this kind of construction~: let $N$ be a normal subgroup of $G$, and let $H=G/N$. Also denote by $U$ the set $G/N$, viewed as a $(G,H)$-biset. As for any $kG$-module $W$, there is an isomorphism of $kH$-modules
$$\Hom_{kG}(kU,W)\cong W^N\mvirg$$
the functor $h_U : \per^+_k(G)\to \per^+_k(H)$ is the fixed points functor $W\mapsto W^N$. \par
Now $W^N\supseteq {\rm tr}_\un^N(W)$, and this inclusion is functorial with respect to $W$~: the functor $W\mapsto {\rm tr}_\un^N(W)$ is a functor from $\gMod{RG}$ to $\gMod{RH}$, which is a subfunctor of the functor $W\mapsto W^N$.\par
A natural question is then to know whether it can happen that the functor $W\mapsto {\rm tr}_\un^N(W)$ preserves direct summands of permutation modules. This is equivalent to requiring that the following condition holds~:
\begin{mth}{Condition}\label{cond}
For any finite $G$-set $X$, the module ${\rm tr}_\un^N(kX)$ isomorphic to a direct summand of a permutation $kH$-module.
\end{mth}
\begin{mth}{Lemma} Let $G$ be a finite group, let $N$ be a normal subgroup of $G$, let $H=G/N$, and let $X$ be a finite $G$-set. Then there is an isomorphism of $kH$-modules
$${\rm tr}_\un^N(kX)\cong\bigoplus_{x\in G\dom X}\Ind_{NG_x/N}^{G/N}(|N_x|k)\mvirg$$ where $G_x$ is the stabilizer of $x$ in $G$, and $N_x=N\cap G_x$.
\end{mth}
\pf The module $(kX)^N$ has a $k$-basis consisting of the elements 
$$b_x=\sum_{n\in N/N_x}nx\mvirg$$ 
for $x\in N\dom X$. This basis is $G/N$-invariant, and the stabilizer of $b_x$ in $G/N$ is equal to $NG_x/N$. Thus
$$(kX)^N\cong \bigoplus_{x\in G\dom X}\Ind_{NG_x/N}^{G/N}k\mpoint$$
The submodule ${\rm tr}_\un^N(kX)$ is generated by the elements ${\rm tr}_\un^N(x)=|N_x|b_x$, for $x\in N\dom X$, and these elements are permuted by $G/N$. The lemma follows.\findemo
In particular, as a $k$-module, 
$${\rm tr}_\un^N(kX)\cong\bigoplus_{x\in N\dom X}(|N_x|k)\mpoint$$
Thus if ${\rm tr}_\un^N(kX)$ is isomorphic to a direct summand of a permutation $kH$-module, since any permutation $kH$-module is free as a $k$-module, it follows that $|N_x|k$ is a projective $k$-module. Hence, if Condition~\ref{cond} holds, then the following condition also holds~:
\begin{mth}{Condition}\label{cond2} For any subgroup $M$ of $N$, the module $|M|k$ is a projective $k$-module.
\end{mth}
Conversely, if Condition~\ref{cond2} holds, and if $X$ is a finite $G$-set, then for any $x\in X$, the $k$-module $|N_x|k$ is projective, hence it is isomorphic to a direct summand of $k$ (since the map $\lambda\mapsto |N_x|\lambda$ is a surjective morphism of $k$-modules from $k$ to $|N_x|k$). It follows that $|N_x|k$ is also isomorphic to a direct summand of $k$ as $k(NG_x/N)$-module, since the action of $NG_x/N$ on $k$ and $|N_x|k$ is trivial. Thus $\Ind_{NG_x/N}^{G/N}(|N_x|k)$ is isomorphic to a direct summand of $\Ind_{NG_x/N}^{G/N}(k)$ as a $k(G/N)$-module. Since this is a permutation module, it follows that ${\rm tr}_\un^N(kX)$ is isomorphic to a direct summand of a permutation $kH$-module. This shows that Conditions~\ref{cond} and~\ref{cond2} are equivalent.
\npar Now saying that for some $m\in \N$, the $k$-module $mk$ is projective, is equivalent to saying that there exists an idempotent $e_m$ of $k$ such that $e_mk$ is equal to the annihilator ${\rm Ann}_k(m)$ of $m$ in $k$. Moreover Condition~\ref{cond2} obviously implies the following~:
\begin{mth}{Condition}\label{cond3} For any prime factor $p$ of $|N|$, the $k$-module $pk$ is projective. Equivalently, there exists an idempotent $e_p$ of $k$ such that ${\rm Ann}_k(p)=e_pk$.
\end{mth}
Conversely, if Condition~\ref{cond3} holds, let $m$ be any integer dividing $|N|$, and denote by $e_m$ the sum $\sum_{p|m}\limits e_p$ of the idempotents $e_p$ corresponding to the {\em distinct} prime factors of $m$. Then $e_m$ is an idempotent~: indeed, if $p$ and $q$ are distinct prime numbers, then $e_pe_q=0$, since it is both a $p$-torsion and a $q$-torsion element. Note that, more generally, if $m$ and $n$ are integers such that $m|n$, then $e_me_n=e_m$. \par
Since moreover 
$$me_m=\sum_{p|m}\frac{m}{p}pe_p=0\mvirg$$
it follows that $e_mk\subseteq {\rm Ann}_k(m)$.\par
Now write $m$ as a product $p_1p_2\cdots p_l$ of (possibly equal) prime numbers. We prove by induction on $l$ that ${\rm Ann}_k(m)=e_mk$. \par
If $l=0$, then $m=1$, and $e_1=0$ generates ${\rm Ann}_k(1)=\zero$. If $l>0$, let $p=p_1$, and $n=m/p$. If $x\in{\rm Ann}_k(m)$, then $pnx=0$, so there exists $y\in k$ such that $nx=e_py$. Thus $e_pnx=nx$, i.e. $n(x-e_px)=0$. By induction hypothesis, there exists $z\in k$ such that $x=e_px+e_nz$. Thus
$$e_mx=e_me_px+e_me_nz=e_px+e_nz=x\mvirg$$
since $e_me_r=e_r$ if $r|m$. It follows that $x\in e_mk$, as was to be shown.\par
This shows that Conditions~\ref{cond},~\ref{cond2}, and~\ref{cond3} are equivalent.

\begin{rem}{Remark} \label{hereditary}These conditions are fulfilled in particular if the ring $k$ is hereditary, since any ideal of $k$ is projective in this case.
\end{rem}

\begin{mth}{Notation} When the normal subgroup $N$ of $G$ fulfills Condition~\ref{cond3}, the functor $W\mapsto{\rm tr}_\un^NW$ from $\per^+_k(G)$ to $\per^+_k(G/N)$ induces an exact functor denoted by $\sigma_{G/N}^G$ from $\comack_k(G/N)$ to $\comack_k(G)$.
\end{mth}
The exactness of $\sigma_{G/N}^G$ follows from the fact that this functor is obtained by {\em pre-composition} with some functor.
\npar\label{definition simple functor} From now on we will assume that $k$ is a field of characteristic $p\geq 0$. Remark~\ref{hereditary} shows that Condition~\ref{cond3} is fulfilled, for any normal subgroup $N$ of $G$.\par
Recall that the simple cohomological Mackey functors for $G$ over $k$ are indexed by pairs $(Q,V)$, where $Q$ is a $p$-subgroup of $G$ (which in the case $p=0$ should be understood as $Q=\un$), and $V$ is a simple $k\sur{N}_G(Q)$-module, where, as usual, $\sur{N}_G(Q)=N_G(Q)/Q$. As an object of $\fun_k^+(G)$, the functor $S_{Q,V}^G$ indexed by the pair $(Q,V)$ can be described as follows~: if $W$ is a direct summand of a permutation $kG$-module, then
$$S_{Q,V}^G(W)=\tr_\un^{\sur{N}_G(Q)}\Hom_k(W[Q],V)\mvirg$$
where $W[Q]$ is the Brauer quotient of $W$ at $Q$. The functorial structure with respect to $W$ is the obvious one. In other words $S_{Q,V}^G(W)$ is the set of $k\sur{N}_G(Q)$-homomorphisms from $W[Q]$ to $V$ which factor through a projective $k\sur{N}_G(Q)$-module.
\begin{mth}{Proposition} \label{sigma simple} Let $k$ be a field. If $N\normal G$, and if $V$ is a simple $k(G/N)$-module,
then
$$\sigma_{G/N}^G(S_{\un,V}^{G/N})=S_{\un,\Inf_{G/N}^GV}^G\mpoint$$
\end{mth}
\pf The value of the functor $S_{\un,V}^{G/N}$ on a direct summand of a permutation $k(G/N)$-module $W$ is equal to
$$S_{\un,V}^{G/N}(W)={\rm tr}_\un^{G/N}\Hom_k(W,V)\mpoint$$
So for $U\in \per^+_k(G)$
\begin{eqnarray*}
\sigma_{G/N}^G(S_{\un,V}^{G/N})(U)&=&S_{\un,V}^{G/N}\big({\rm tr}_\un^N(U)\big)\\
&=&{\rm tr}_\un^{G/N}\Hom_k\big({\rm tr}_\un^N(U),V\big)\mpoint
\end{eqnarray*}
On the other hand
$$S_{\un,\Inf_{G/N}^GV}^G(U)={\rm tr}_\un^G\Hom_k(U,\Inf_{G/N}^GV)\mpoint$$
Let $i: \tr_\un^NU\hookrightarrow U$ denote the inclusion map, and let $s:U\to \tr_\un^NU$ be a $k$-linear map such that $s\circ i=\Id$. If $\psi\in\Hom_k(\tr_\un^NU,V)$, then $\psi=\widetilde{\psi}\circ i$, where $\widetilde{\psi}=\psi\circ s\in\Hom_k(U,\Inf_{G/N}^GV)$. Setting $\theta=\tr_\un^{G/N}\psi$, for any $u\in U$
\begin{eqnarray*}
\theta\circ\tr_\un^N(u)&=&\sum_{g\in G/N}\sum_{n\in N}g\widetilde{\psi}(g^{-1}nu)\\
&=&\sum_{g\in G/N}\sum_{n\in N}gn^{-1}\widetilde{\psi}(g^{-1}nu)\\
&=&\sum_{g\in G}g\psi(g^{-1}u)\\
&=&\tr_\un^G(\widetilde{\psi})(u)\mvirg
\end{eqnarray*}
thus $\theta\circ\tr_\un^N=\tr_\un^G(\widetilde{\psi)}$. This shows that the map
$$I:\Hom_k(\tr_\un^NU,V)\to \Hom_k(U,\Inf_{G/N}^GV)$$
sending $\theta$ to $\theta\circ \tr_\un^N$ maps $\tr_\un^{G/N}\Hom_k(\tr_\un^NU,V)$ into $\tr_\un^G\Hom_k(U,\Inf_{G/N}^GV)$. The map $I$ is obviously injective, as $\tr_\un^N:U\to \tr_\un^NU$ is surjective.\par
Conversely, let $\varphi\in \Hom_k(U,\Inf_{G/N}^GV)$, and set $\Xi=\tr_\un^G\varphi$. Then for $u\in U$
\begin{eqnarray*}
\Xi(u)&=&\sum_{g\in G}g\varphi(g^{-1}u)\\
&=&\sum_{g\in G/N}\sum_{n\in N}gn\varphi(n^{-1}g^{-1}u)\\
&=&\sum_{g\in G/N}\sum_{n\in N}g\varphi(g^{-1}n^{-1}u)\\
& &\hbox{(since }n^{-1}g^{-1}w=g^{-1}{\cdot}gn^{-1}g^{-1}\cdot u\hbox{)}\\
&=&\big({\rm tr}_\un^{G/N}(\varphi\circ i)\big)(\tr_\un^Nu)\mvirg
\end{eqnarray*}
so $\Xi=I\big({\rm tr}_\un^{G/N}(\varphi\circ i)\big)$, showing that $I$ induces an isomorphism
$$\tr_\un^{G/N}\Hom_k(\tr_\un^NU,V)\to\tr_\un^G\Hom_k(U,\Inf_{G/N}^GV)\mvirg$$
which is obviously functorial in $U$. This completes the proof.
\findemo
\begin{mth}{Proposition} \label{consequence sigma}Let $G=N\rtimes H$ be the semidirect product of $N$ by a group~$H$.
\begin{enumerate}
\item If $V$ is a $kH$-module, let $\tilde{V}$ be the $kG$-module $\Inf_{G/N}^G\Iso_H^{G/N}V$. Then the restriction of $\tilde{V}$ to $H$ is isomorphic to $V$.
\item The composition $\Res_H^G\sigma_{G/N}^G\Iso_{H}^{G/N}$ is isomorphic to the identity functor of $\comack_k(H)$.
\item Let $V$ and $W$ be simple $kH$-modules. Then, for any $n\in \N$, the restriction from $G$ to $H$ induces a split surjection
$$r_H^G:\Ext^n_{\comack_k(G)}(S_{\un,\tilde{V}}^G,S_{\un,\tilde{W}}^G)\to\Ext^n_{\comack_k(H)}(S_{\un,V}^H,S_{\un,W}^H)\mpoint$$

\end{enumerate}
\end{mth}
\pf Assertion~1 is obvious. For Assertion~2, let $W$ be a direct summand of a permutation $kH$-module, and let $W'$ be the $kH$-module defined by
$$W'=\Iso_{G/N}^H{\rm tr}_\un^N\Ind_H^GW\mpoint$$
The set $N$ is a set of representatives of $G/H$, so
$$\Ind_H^GW\cong \bigoplus_{n\in N}n\otimes W\mpoint$$
Moreover, for any $n\in N$ and $w\in W$,
$${\rm tr}_\un^N(n\otimes w)=\sum_{x\in N}xn\otimes w=\sum_{x\in N}x\otimes w\mpoint$$
This shows that the map $\theta:w\in W\mapsto \dsp\sum_{x\in N}\limits x\otimes w\in W'$ is a $k$-linear isomorphism from $W$ to $W'$. Moreover, for any $h\in H$
\begin{eqnarray*}
h(\sum_{x\in N}x\otimes w)&=&\sum_{x\in N}hx\otimes w\\
&=&\sum_{x\in N}hxh^{-1}\cdot h\otimes w=\sum_{x\in N}xh\otimes w\\
&=&\sum_{x\in N}x\otimes hw\mpoint
\end{eqnarray*}
Thus $\theta$ is actually an isomorphism of $kH$-modules $W\to W'$, which is obviously functorial with respect to $W$. So the functor $F=\Iso_{G/N}^H{\rm tr}_\un^N\Ind_H^G$ is isomorphic to the identity functor of $\per^+_k(H)$. Assertion~2 follows, since the functor $\Res_H^G\sigma_{G/N}^G\Iso_{H}^{G/N}$ is the endofunctor of $\comack_k(H)$ obtained by composition with $F$.\par
Assertion~1 implies that $\Res_H^GS_{\un,\tilde{V}}^G\cong S_{\un,V}^H$, so the restriction functor $\Res_H^G$ induces a $k$-linear map
$$r_H^G:\Ext^n_{\comack_k(G)}(S_{\un,\tilde{V}}^G,S_{\un,\tilde{W}}^G)\to\Ext^n_{\comack_k(H)}(S_{\un,V}^H,S_{\un,W}^H)\mpoint$$
Conversely, the functor $\sigma_{G/N}^G\Iso_{H}^{G/N}$ is an exact functor from $\comack_k(H)$ to $\comack_k(G)$, which sends $S_{\un,V}^H$ to $S_{\un,\tilde{V}}^G$, by Proposition~\ref{sigma simple}. This yields a $k$-linear map
$$s_H^G:\Ext^n_{\comack_k(H)}(S_{\un,V}^H,S_{\un,W}^H)\to \Ext^n_{\comack_k(G)}(S_{\un,\tilde{V}}^G,S_{\un,\tilde{W}}^G)\mpoint$$
Since by Assertion~2, the functor $\Res_H^G\sigma_{G/N}^G\Iso_{H}^{G/N}$ is isomorphic to the identity functor of $\comack_k(H)$, the composition $r_H^G\circ s_H^G$ is an isomorphism, so $r_H^G$ is split surjective.\findemo

%%%%%%%%%%%%%%%%%%%%%%%%%%%%%%%%%%%%%%%%%%%%%%%%%%%%%%%%
\section{Extensions of cohomological Mackey functors}%%%
%%%%%%%%%%%%%%%%%%%%%%%%%%%%%%%%%%%%%%%%%%%%%%%%%%%%%%%%
\begin{mth}{Hypothesis}\label{hyp3.1} From now on, we assume that $k$ is a field of positive characteristic $p$ and that $G$ is a finite $p$-group.
\end{mth}
In this section we recall some notation and results from \cite{cohocplx} about extensions of cohomological Mackey functors for elementary abelian $p$-groups. When $H$ is a subgroup of $G$, we denote by $K_G(H)$ the set of complements of $H$ in~$G$, i.e. the set of subgroups $T$ of $G$ such that $H\oplus T=G$.
\npar (\cite{cohocplx} Corollary 8.3) When $X$ is a subgroup of order $p$ of $G$, there exist a unique cohomological Mackey functor $\bin{X}{\un}{G}$, up to isomorphism, which fits into a non split exact sequence in $\comack_k(G)$ of the form
\begin{equation}\label{ext1SXS1}0\to S_\un^G\to \bin{X}{\un}{G}\to S_X^G\to 0\mpoint\end{equation}
There is also a unique cohomological Mackey functor $\bin{\un}{X}{G}$, up to isomorphism, which fits into a non split exact sequence in $\comack_k(G)$ of the form
\begin{equation}\label{ext1S1SX}0\to S_X^G\to \bin{\un}{X}{G}\to S_\un^G\to 0\mvirg\end{equation}
and $\bin{\un}{X}{G}$ is isomorphic to the dual $\bin{X}{\un}{G}^*$ of $\bin{X}{\un}{G}$.\par
We denote by $\gamma_X^G\in\Ext^2_{\comack_k(G)}(S_\un^G,S_\un^G)$ the class of the exact sequence 
$$\Gamma_X^G: 0\to S_\un^G\to\bin{X}{\un}{G}\to\bin{\un}{X}{G}\to S_\un^G\to 0$$
obtained by splicing the two previous short exact sequences (\cite{cohocplx} Notation~7.4).
\npar \label{tauphi}(\cite{cohocplx} Section~14) When $p>2$, and $\varphi\in\Hom(G,k^+)$, let $U_\varphi^G$ denote the $kG$-module $k\oplus k$, on which $G$ acts by
$$\forall g\in G,\;\forall (x,y)\in k^2,\;\;g(x,y)=\big(x+\varphi(g)y,y\big)\mpoint$$
Let $T_\varphi^G$ denote the unique Mackey functor for $G$ over $k$ such that $T_\varphi^G(H)$ is zero if $H$ is a non trivial subgroup of $G$, and $T_\varphi^G(\un)\cong U_\varphi^G$. The functor $T_\varphi^G$ is cohomological, and fits in an exact sequence
$$0\to S_\un^G\to T_\varphi^G\to S_\un^G\to 0\mvirg$$
in $\comack_k(G)$, which is non split if $\varphi\neq 0$. We denote by $\tau_\varphi^G$ the class of this extension in $\Ext^1_{\comack_k(G)}(S_\un^G,S_\un^G)$.
 When $\varphi\in\Hom(G,\F_p^+)$, we denote by the same symbol the composition of $\varphi$ with the inclusion $\F_p^+\hookrightarrow k^+$, and by $\tau_\varphi$ the corresponding element of $\Ext^1_{\comack_k(G)}(S_\un^G,S_\un^G)$.
\npar The following conjecture was proposed in \cite{cohocplx}, and proved there for $p=3$ (\cite{cohocplx}~Theorem~14.2)~:
\begin{mth}{Conjecture} \label{conjecture}Let $k$ be a field of odd characteristic $p$, and $G\cong (C_p)^r$. Then~:
\begin{enumerate}
\item The algebra $\Ee=\Ext^*_{\comack_k(G)}(S_\un^G,S_\un^G)$ is generated by the elements $\tau_\varphi^G$ in degree 1, for $\varphi\in \Hom(G,k^+)$, and by the elements $\gamma_X^G$ in degree~2, for $X\leq G$ with $|X|=p$.
\item The Poincar\'e series for $\Ee$ is equal to
$$\frac{1}{(1-t)\big(1-t-(p\!-\!1)t^2\big)\big(1-t-(p^2\!-\!1)t^2\big)\ldots\big(1-t-(p^{r-1}\!-\!1)t^2\big)}\mpoint$$
\end{enumerate}
\end{mth}
\begin{mth}{Proposition}{\rm[\cite{cohocplx} Proposition 8.7 and Proposition 10.1]} \label{index p elemab}Let $k$ be a field of characteristic~$p>0$, let $G$ be an elementary abelian $p$-group, and let $H$ be a subgroup of index $p$ in $G$. Set $I=\Ind_H^GS_\un^H$.
\begin{enumerate}
\item Let $R$ and $S$ denote respectively the radical and the socle of $I$, as an object of $\comack_k(G)$. Then $I\supset R\supseteq S\supset \zero$, and $I/R\cong S\cong S_\un^G$. Moreover
$$R/S\cong L\oplus\dirsum{X\in K_G(H)}S_X^G\mvirg$$
where $L$ is a functor all of whose composition factors are isomorphic to $S_\un^G$, with multiplicity $p-2$.
\item Let $Y\in K_G(H)$. The functor $R$ has a subfunctor~$J$ isomorphic to $\iota_{G/Y}^G(S_\un^{G/Y})$, and there is an isomorphism
$$R/J\cong L\oplus \dirsum{X\in K_G(H)-\{Y\}}S_X^G\mvirg$$
\end{enumerate}
\end{mth}
\begin{mth}{Corollary} {\rm[\cite{cohocplx} Corollary 10.3]}\label{long exact sequence} With the same notation, there is a long exact sequence of extension groups
$$\xymatrix{
\cdots\ar[r]&*!U(0.6){L(n-1)\oplus\dirsum{X\in\mathcal{X}}E_G(n-2)}\ar[r]&E_G(n)\ar[r]&E_H(n)\ar`r[d]`[l]`[dlll]`[dll][dll]\\
&*!U(0.6){L(n)\oplus\dirsum{X\in \mathcal{X}}E_G(n-1)}\ar[r]&E_G(n+1)\ar[r]&E_H(n+1) \cdots\mvirg\\
}
$$
where $E_G(n)=\Ext^n_{\comack_k(G)}(S_\un^G,S_\un^G)$, $E_H(n)=\Ext^n_{\comack_k(H)}(S_\un^H,S_\un^H)$, $L(n)=\Ext^n_{\comack_k(G)}(L,S_\un^G)$, and $\mathcal{X}=K_G(H)-\{Y\}$.
\end{mth}
Moreover, it is easy to check that the map $E_G(n)\to E_H(n)$ in this corollary is induced by the restriction functor $\Res_H^G$, since $\Res_H^GS_\un^G\cong S_\un^H$. As $k$ is a field by Hypothesis~\ref{hyp3.1}, the conclusion of Proposition~\ref{consequence sigma} holds and this map is (split) surjective. Thus~:
\begin{mth}{Corollary} \label{short exact sequence}With the same notation, for any $n\in\N$, there is a short exact sequence of extension groups
$$0\longrightarrow L(n-1)\oplus\mathop{\oplus}_{X\in\mathcal{X}}\limits E_G(n-2)\to E_G(n)\to E_H(n)\to 0\mpoint$$
\end{mth}

%%%%%%%%%%%%%%%%%%%%%%%%%%%%%%%%%%%%%%%%%%%%%%%%%%%%%%%%%%%%%%%%%%%%%%%%%
\section{Mackey functors concentrated at $\un$ and relations in~$\Ee$}%%%
%%%%%%%%%%%%%%%%%%%%%%%%%%%%%%%%%%%%%%%%%%%%%%%%%%%%%%%%%%%%%%%%%%%%%%%%%
We first consider some generalizations of the functor $L$ of Proposition~\ref{index p elemab}~:
\begin{mth}{Definition} Let $k$ be a field of characteristic $p>0$, and $G$ be a finite $p$-group. 
\begin{itemize}
\item A Mackey functor $M$ for $G$ over $k$ is said to be {\em concentrated at $\un$} if $M(H)=\zero$ for any non trivial subgroup $H$ of $G$.
\item A $kG$-module $V$ is said to {\em have zero traces} if  
$$\tr_\un^X(V)=\zero$$
for any non trivial subgroup $X$ of $G$.
\end{itemize}
\end{mth}
\begin{rem}{Remark} 0) A Mackey functor concentrated at \un~is cohomological.\par
1) A (finitely generated) Mackey functor for $G$ over $k$ is concentrated at $\un$ if and only if all its composition factors are isomorphic to~$S_\un^G$. \par
2) By transitivity of traces, a $kG$-module $V$ has zero traces if and only if $\tr_\un^X(V)=\zero$ for any subgroup $X$ {\em of order $p$} of $G$.
\end{rem}

\begin{mth}{Proposition} \label{concentrated}Let $k$ be a field of characteristic $p$, and $G$ be a finite $p$-group.
\begin{enumerate}
\item  Let $M$ be a Mackey functor for $G$ over $k$. If $M$ is concentrated at~$\un$, then $M$ is cohomological, and the $kG$-module $V=M(\un)$ has zero traces.
\item If $V$ is a $kG$-module having zero traces, then there is a unique Mackey functor $\widehat{V}$ for $G$ over $k$ such that $\widehat{V}$ is concentrated at $\un$ and $\widehat{V}(\un)\cong V$ as $kG$-modules.
\item If $M$ is an object of $\comack_k(G)$, let $M^0$ denote the $kG$-submodule
$$M^0=\bigcap_{1<X\leq G}\Ker\;t_\un^X$$
of $M(1)$. Then $M^0$ has zero traces, and $\widehat{M^0}$ is the largest subfunctor of $M$ concentrated at $\un$.
\item The correspondences $M\mapsto M(\un)$ and $V\mapsto\widehat{V}$ are mutual inverse equivalences of categories between the full subcategory of $\comack_k(G)$ whose objects are concentrated at $\un$, and the full subcategory of $\gMod{kG}$ whose objects are modules with zero traces.
\end{enumerate}
\end{mth}
\pf For Assertion 1, observe that for any subgroup $X$ of $G$, and any $v\in V=M(\un)$
$$r_\un^Xt_\un^Xv=\sum_{x\in X}x\cdot v=\tr_\un^X(v)\mpoint$$
This has to be zero if $X$ is non-trivial, since $M(X)=\zero$ by assumption. Moreover $M$ is obviously cohomological.\par
For Assertion 2, it is straightforward to check that if $V$ has zero traces, then the assignments $M(\un)=V$ and $M(H)=\zero$ for $\un<H\leq G$, define a Mackey functor $M$ for~$G$ over~$k$. This proves the existence part of Assertion~2. Uniqueness is straightforward. \par
Assertion 3 is also straightforward, and Assertion 4 follows easily.
\findemo
\begin{rem}{Example} \label{example zero traces}Let $H$ be a subgroup of index $p$ of $G$, and let $W=\Ind_H^Gk$. Let $\varepsilon : W\to k$ denote the augmentation map, and set $V=\Ker\,\varepsilon$. Then the $kG$-module $V$ has zero traces~: indeed, if $X$ is a subgroup of order $p$ of $G$, then either $X\leq H$, and then $\tr_\un^X(W)=\zero$, since $X$ acts trivially on $W$. Or $X\nleq H$, and then $\Res_X^GW\cong\Ind_\un^Xk\cong kX$, so there is an exact sequence
$$0\to \Res_X^GV\to kX\to k\to 0\mvirg$$
showing that $\tr_\un^X(V)=\zero$ also in this case.\par
The module $\Ind_H^Gk$ is inflated from the free module of rank 1 for the cyclic group $G/H$. Hence, it is uniserial, and all its subquotients are indecomposable, and characterized by their dimension, up to isomorphism. We denote by $U_a$ the subquotient of dimension $a$, for $a\in\{1,\ldots p\}$. If $a\leq p-1$, the module $U_a$ is isomorphic to a submodule of $V$, hence it has zero traces, and we denote by $T_a$ the functor $\widehat{U_a}$. \par
Let $\varphi\in\Hom(G,k^+)$ with kernel $H$. Then the functor $T_\varphi^G$ introduced in Paragraph~\ref{tauphi} is isomorphic to $T_2$. Similarly, the functor $L$ of Proposition~\ref{index p elemab} is isomorphic to $T_{p-2}$, since $L(1)\cong U_{p-2}$. We also denote by $M$ the functor $T_{p-1}=\widehat{V}$. The two short exact sequences 
\begin{equation}\label{kUaUa+1}
0\to k\to U_{a+1}\to U_a\to 0\;\;,\;\;\hbox{and}\;\;0\to U_a\to U_{a+1}\to k\to 0
\end{equation}
of $kG$-modules yield corresponding exact sequences
\begin{equation}
0\to S_\un^G\to T_{a+1}\to T_a\to 0\;\;,\;\;\hbox{and}\;\;0\to T_a\to T_{a+1}\to S_\un^G\to 0
\end{equation}\label{kTaTa+1}
of cohomological Mackey functors for $G$ over $k$. In particular,
\begin{equation}\label{kVU}
0\to k\to V\to U\to 0\;\;,\;\;\hbox{and}\;\;0\to U\to V\to k\to 0
\end{equation}
yield
\begin{equation}\label{kML}
0\to S_\un^G\to M\to L\to 0\;\;,\;\;\hbox{and}\;\;0\to L\to M\to S_\un^G\to 0\,.
\end{equation}
%\begin{rem}{Example} \label{define M}With the notation of Proposition~\ref{index p elemab}, set $M=\widehat{I^0}$. Then $M$ is a proper subfunctor of $I$, so $M\subseteq R$, and $M$ is non zero, so $M\supseteq S$, and clearly $M/S\cong L$. In other words $M(\un)\leq \Ind_H^Gk$ is the kernel $V$ of the augmentation map $\Ind_H^Gk\to k$, like in Example~\ref{example zero traces}. By Proposition~\ref{concentrated}, the first exact sequence in~\ref{kVU} corresponds to the short exact sequence
%$$0\to S_\un^G\to M\to L\to 0\mvirg$$
%and the second one to an exact sequence
%$$0\to L\to M\to S_\un^G\to 0\mvirg$$
%so $L$ is also isomorphic to a subfunctor of $M$ (actually the radical of $M$).
We will denote by $\tau_{\varphi,L}$ the element of $\Ext^1_{\comack_k(G)}(L,S_\un^G)$, respectively $\tau_{L,\varphi}$ the element of $\Ext^1_{\comack_k(G)}(S_\un^G,L)$ corresponding to this last two short exact sequences.\par
Assertion~2 of Proposition~\ref{index p elemab} can now be rephrased as follows~: there exists a subfunctor $\Sigma$ of $R$, containing $S$, such that
$$R=M+\Sigma,\;\;S=M\cap\Sigma,\;\;\Sigma/S\cong \dirsum{X\in K_G(H)}S_X^G\mpoint$$
Setting $K=L+\Sigma$, this gives the following diagram of subfunctors of $I$~:
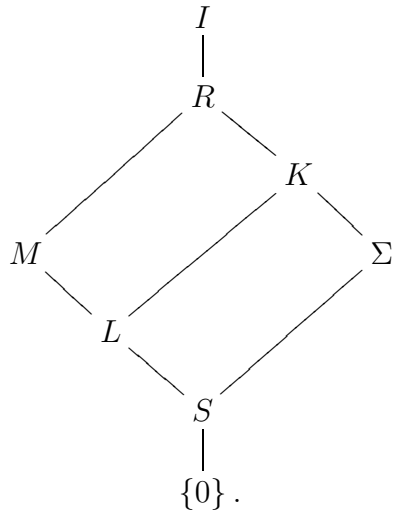
$$\xymatrix@C=3ex@R=3ex{
&&I\ar@{-}[d]&&\\
&&R\ar@{-}[ddll]\ar@{-}[dr]&&\\
&&&K\ar@{-}[ddll]\ar@{-}[dr]&\\
M\ar@{-}[dr]&&&&\Sigma\ar@{-}[ddll]\\
&L\ar@{-}[dr]&&&\\
&&S\ar@{-}[d]&&\\
&&\zero\makebox[0pt]{\ \ .}&&
}
$$
\end{rem}
We give now a series of relations between elements in $\Ext^1_{\comack_k(G)}(S_\un^G,S_\un^G)$ and $\Ext^2_{\comack_k(G)}(S_\un^G,S_\un^G)$. For simplicity, we will drop the exponent $G$, and write $\tau_\varphi$, $\gamma_X$, $S_\un$, $T_\varphi$, $U_\varphi$ instead of $\tau_\varphi^G$, $\gamma_X^G$, $S_\un^G$, $T_\varphi^G$, $U_\varphi^G$ respectively. We start with a linearity relation between the $\tau$'s

\begin{mth}{Lemma}\label{TauLinRel}
Let $\varphi$ and $\psi$ be two morphisms in $\Hom(G,k^+)$. Then $\tau_\varphi+\tau_\psi=\tau_{\varphi+\psi}$ in $\Ext^1_{\comack_k(G)}(S_\un,S_\un)$.
\end{mth}
\pf
Let $0\to S_\un\to T_\varphi\to S_\un\to 0$ be the representative for $\tau_\varphi$ and $0\to S_\un\to T_\psi\to S_\un\to 0$ be the representative for $\tau_\psi$ described in Paragraph~\ref{tauphi}. Given that both $\tau_\varphi$ and $\tau_\psi$ are concentrated at the trivial subgroup, this will also be the case for $\tau_\varphi+\tau_\psi$. \par
We construct a representative $0\to S_\un\to T\to S_\un\to 0$ for the sum as in Lemma~\ref{CompExt1} working with the $G$-modules $U_\varphi$ and $U_\psi$ that give the functors $T_\varphi$, respectively $T_\psi$. Then $T(1)$ is given by the following sequence of pushout and pullback
$$\xymatrix{
k\oplus k\ar[r]^(0.4){i\oplus j}\ar[d]_{\pi_1+\pi_2}&U_\varphi\oplus U_\psi\ar[d]\\
k\ar[r]&U'}\,,\qquad\xymatrix{
U'\ar[r]^(0.4){\widetilde{s\oplus t}}&k\oplus k\\
T(1)\ar[r]\ar[u]&k\ar[u]_{\Delta}}\,.
$$
The same computation as in Lemma~\ref{CompExt1} gives 
$$T(1)\simeq \{(a_1,a_2,b_1,b_2,c)\in k^5|a_2=b_2\}/\{(0,d_1,0,d_2,-d_1-d_2)|d_1,d_2\in k\}\,.$$
This module has dimension $2$, the action of $G$ on the class of $(a_1,a_2,b_1,b_2,c)$ is induced by the action on $U_\varphi$ for the first two terms, the action on $U_\psi$ for the next two terms and is trivial on the last. More precisely, $$g[(a_1,a_2,b_1,b_2,c)]=[(a_1+\varphi(g)a_2,a_2,b_1+\psi(g)b_2,b_2,c)]\,.$$ 
Thus $T(1)$ is a $G$-module with $k$-basis $\{[(0,1,0,1,0)],[(0,0,0,0,1)]\}$ and action given by 
$$g[(0,0,0,0,1)]=[(0,0,0,0,1)]$$ 
and 
$$g[(0,1,0,1,0)]\!=\![(\varphi(g)1,1,\psi(g),1,0)]\!=\![(0,1,0,1,0)]+[(0,0,0,0,\varphi(g)+\psi(g))]\,.$$ 
This means that $T(1)\cong U_{\varphi+\psi}$ and that $T\cong T_{\varphi+\psi}$.
\findemo
The following lemmas give relations involving $(\tau_\varphi)^2$. We obtain different relations in the cases $p=3$ and $p\ge 5$. This is tidily related to the existence of the functor $T_3$ only when $p\ge 5$. This functor is constructed in Example \ref{example zero traces} as an application of Proposition~\ref{concentrated}~. The main ingredient used in both lemmas is Lemma~12.2 of~\cite{cohocplx}, that we use to detect zero elements in $\Ext^2_{\comack_k(G)}(S_\un,S_\un)$.
\begin{mth}{Lemma}\label{TauSqIsZero} Let $k$ be a field of characteristic $p\geq 5$ and let $\varphi\in\Hom(G,k^+)$. Then $(\tau_\varphi)^2=0$ in $\Ext^2_{\comack_k(G)}(S_\un,S_\un)$.
\end{mth}
\pf Recall from Example~\ref{example zero traces} that $T_\varphi=\bn{S_\un}{S_\un}\cong T_2$, when $H=\Ker\,\varphi$. Since $p\geq 5$, we have $p-1\geq 3$, so the module $U_3$ defined in Example~\ref{example zero traces} has zero traces. This module has a filtration
$$\zero\subset k\subset U_2\subset U_3\mvirg$$
such that $U_3/k\cong U_2$ and $U_3/U_2\cong k$. It follows that $T_3$ has a filtration
$$\zero\subset S_\un\subset T_2\subset T_3$$
such that $T_3/S_\un\cong T_2$ and $T_3/T_2\cong S_\un$. 
Applying Lemma~12.2 of~\cite{cohocplx} to this filtration shows that the sequence
$$0\to S_\un\to \bn{S_\un}{S_\un}\to\bn{S_\un}{S_\un}\to S_\un\to 0$$
represents the zero element of $\Ext^2_{\comack_k(G)}(S_\un,S_\un)$. But this represents precisely~$(\tau_\varphi)^2$.~\findemo
In the case $p=3$, using the notation of Example \ref{example zero traces}, we have $M=T_2$ and $L=S_\un$. Hence, the construction in the proof of previous lemma cannot be used. We use instead the decomposition of the functor $I$ from Proposition~\ref{index p elemab}.
\begin{mth}{Hypothesis} We assume from now on that $G\cong (C_p)^r$ is an elementary abelian $p$-group of rank $r$.
\end{mth}
\begin{mth}{Lemma}\label{TauSqForp3} Let $k$ be a field of characteristic $p=3$. Let $G$ be an elementary abelian $p$-group, and $\varphi\in\Hom(G,\F_p^+)$. Then 
$$(\tau_\varphi)^2=-\sum_{X\nleq\Ker\,\varphi}\limits\hat\gamma_X$$ 
in $\Ext^2_{\comack_k(G)}(S_\un,S_\un)$.
\end{mth}
\pf
If $\varphi=0$, there is nothing to prove, since $\tau_\varphi=0$ and the summation in the right hand side is empty. So we assume $\varphi\neq 0$, and set $H=\Ker\varphi$. With the notation of Proposition \ref{index p elemab}, let $I\supset R\supset S_\un$ be a filtration of the functor $I=\Ind_H^G S_\un^H$. Then we have $I/R\simeq S_\un$ and, using that $R/S_\un\cong S_\un\oplus\dirsum{X\in K_G(H)}S_X$, we obtain $I/S_\un\simeq \bn{S_\un}{S_\un\oplus\dirsum{X\in K_G(H)}S_X}$. Moreover, we have that $R\cong\bn{S_\un\oplus\dirsum{X\in K_G(H)}S_X}{S_\un}$. Lemma~12.2 of~\cite{cohocplx} applied to this filtration gives that the exact sequence
$$0\to S_\un\to \bn{S_\un\oplus\dirsum{X\in K_G(H)}S_X}{S_\un}\to\bn{S_\un}{S_\un\oplus\dirsum{X\in K_G(H)}S_X}\to S_\un\to 0$$
represents the zero element of $\Ext^2_{\comack_k(G)}(S_\un,S_\un)$. The result follows as this sequence also represents $(\tau_\varphi)^2+\sum_{X\in K_G(H)}\limits\hat\gamma_X$.
\findemo
We prove now two commutation relations satisfied by the elements $\tau_\varphi$  and $\gamma_X$ in $\Ext^3_{\comack_k(G)}(S_\un,S_\un)$~:
\begin{mth}{Lemma}\label{GammaCommTau}
Let $k$ be a field of characteristic $p\ge 3$ and $G$ be an elementary abelian $p$-group. Then
$$\gamma_X\tau_\varphi=\tau_\varphi\gamma_X$$
for all $\varphi:G\to k^+$ and $X\le\Ker\,\varphi$, $|X|=p$.
\end{mth}
\pf We fix $\varphi:G\to k^+$ and $X\le\Ker\,\varphi$, $|X|=p$.
Let $M$ be the functor defined by the following data. First set $M(1)=U_\varphi$. Recall that this means $M(1)\cong k\oplus k$ as $k$-vector spaces and the action by $g\in G$ is given by $g(x,y)=(x+\varphi(g)y,y)$. Then set $M(X)=k\oplus k=\{(s,t)|s,t\in k\}$ with trivial $G$-action and $M(H)=0$ for all subgroups $H$ of $G$ different from $\un$ and $X$. Also set $t_\un^X(1,0)=(0,0)$, $t_\un^X(0,1)=(1,0)$, $r_\un^X(1,0)=(0,0)$ and $r_\un^X(0,1)=(1,0)$. It is easy to check that $M$ is a cohomological Mackey functor. The fact that $X\le\Ker\,\varphi$ gives that the action of $u\in X$ is trivial on $U_\varphi$ and, thus $\sum_{u\in X}\limits u\cdot(x,y)=p\cdot(x,y)=0=r_\un^Xt_\un^X(x,y)$. \par
The functor $M$ has a subfunctor isomorphic to $\bn{S_X}{S_\un}$ and the quotient by this subfunctor is isomorphic to $\bn{S_\un}{S_X}$. 
Define the Mackey functor $M'$ by $M'(1)=U_\varphi$, $M'(X)=k$, $r_\un^X(1)=(1,0)$ and $t_\un^X(x,y)=0$. It is straight forward that $M'$ is a cohomological Mackey functor, that $T_\varphi$ is a subfunctor of $M'$ and that $M'=\begin{pmatrix}{\bn{S_X}{S_\un}}\\{S_\un}\end{pmatrix}$. Also, define the Mackey functor $M''$ by $M''(1)=U_\varphi$, $M''(X)=k$, $r_\un^X(1)=(0,0)$ and $t_\un^X(x,y)=y$. Again, it is straight forward that $M''$ is a cohomological Mackey functor, that $T_\varphi$ is a factor of $M''$ when factoring out the socle $S_X$ and that $M''=\begin{pmatrix}{S_\un}\\{\bn{S_\un}{S_X}}\end{pmatrix}$. Moreover, we get an exact sequence
$$0\to S_\un\to M'\to M\to M''\to S_\un\to 0\,.$$ 
It is now straightforward to show that this sequence is Yoneda equivalent to 
$$0\to S_\un\to\bn{S_X}{S_\un}\to\bn{S_\un}{S_X}\to\bn{S_\un}{S_\un}\to S_\un\to 0$$
via the canonical projections $M'\to\bn{S_X}{S_\un}$, $M\to\bn{S_\un}{S_X}$ and, respectively, $M''\to\bn{S_\un}{S_\un}$ and also equivalent to
$$0\to S_\un\to\bn{S_\un}{S_\un}\to\bn{S_X}{S_\un}\to\bn{S_\un}{S_X}\to S_\un\to 0$$
via the canonical injections $\bn{S_\un}{S_\un}\to M'$, $\bn{S_X}{S_\un}\to M$ and, respectively, $\bn{S_\un}{S_X}\to M''$. Given that the last two exact sequences represent $\gamma_X\tau_\varphi$, respectively, $\tau_\varphi\gamma_X$, the result follows. \findemo
\begin{mth}{Lemma}\label{SumGammaCommTau}
Let $k$ be a field of characteristic $p\ge 3$ and $G$ be an elementary abelian $p$-group. Then
$$(\sum_{X\nleq\Ker\,\varphi}\gamma_X)\tau_\varphi =\tau_\varphi(\sum_{X\nleq\Ker\,\varphi} \gamma_X)\mpoint$$
\end{mth}
\pf Again if $\varphi=0$, there is nothing to prove, so we assume $\varphi\neq 0$, and set $H=\Ker\,\varphi$. Let $T=T_\varphi^G$ be the functor described in Paragraph \ref{tauphi} and let $L$, respectively $M$, be the functors described in Example~\ref{example zero traces}.
\sp
{\bf $\bullet$ First step}.
$$(*)\quad \tau_{L,\varphi} \tau_\varphi + \sum_{X\in K_G(H)} \gamma_{L,X} = 0$$
where the notation is as follows~: the functor $L$ was defined in Proposition~\ref{index p elemab}, and the element $\tau_{L,\varphi}$ in Example~\ref{example zero traces}. \par
Let $i: S_\un\to L$ be the inclusion map corresponding to the isomorphism from $S_\un$ to the socle of $L$. Taking the image under $i$ of Extension~\ref{ext1SXS1} yields the following diagram
$$\xymatrix{
0\ar[r]&S_\un\ar[r]\ar[d]^-i&\bn{S_X}{S_\un}\ar[r]\ar[d]&S_X\ar[r]\ar@{=}[d]&0\\
0\ar[r]&L\ar[r]&\bn{S_X}{L}\ar[r]&S_X\ar[r]&0\;\;\\
}
$$
whose bottom line can be spliced with Extension~\ref{ext1S1SX} to give the exact sequence
$$0\to L\to\bn{S_X}{L}\to\bn{S_\un}{S_X}\to S_\un\to 0\;\;,$$
which defines the element $\gamma_{L,X}\in\Ext^2_{\comack_k(G)}(S_\un,L)$.\par
We apply Lemma 12.2 in \cite{cohocplx} to the sequence of functors
$$\zero\subset L \subset R \subset I\;\;,$$
where $R$ and $I$ are defined in Proposition~\ref{index p elemab}. The quotient $R/L$ is isomorphic to $S_\un\oplus\dirsum{X\in K_G(H)}S_X$, so $R$ can be represented by $\bn{S_\un\oplus\dirsum{X\in K_G(H)}S_X}{L}$.

The quotient $I/R$ is isomorphic to $S_\un$. The quotient $I/L$ is isomorphic to $\displaystyle\bn{S_\un}{S_\un\oplus\mathop{\oplus}_{X\in K_G(H)}\limits S_X}$.
Thus one gets a four terms exact sequence  
$$0\to L\to\bn{S_\un\oplus\dirsum{X\in K_G(H)}S_X}{L}\to\bn{S_\un}{S_\un\oplus\dirsum{X\in K_G(H)}S_X}\to  S_\un\to 0$$
representing a zero extension in $\Ext^2(S_\un,L)$, which is in fact (*).
\sp
{\bf $\bullet$ Second step}.
$$(**)\quad \tau_{\varphi,L}\gamma_{L,X}=\tau_\varphi\gamma_X$$
where $\tau_{\varphi,L}=(0\to S_\un\to M\to L\to 0)$

The left hand side of (**) is represented by 
$$0\to S_\un \to  M \to  \bn{S_X}{L} \to  \bn{S_\un}{S_X} \to  S_\un\to 0\;\;,$$ 
and the right hand side is represented by 
$$0\to S_\un \to  T \to  \bn{S_X}{S_\un} \to  \bn{S_\un}{S_X} \to  S_\un\to 0\;\;.$$ 
\par
One constructs an equivalence of these extensions by taking  inclusions maps $T\hookrightarrow M$, respectively $\bn{S_X}{S_\un}\hookrightarrow \bn{S_X}{L}$, and the identity elsewhere. The equality in (**) follows.
\sp
{\bf $\bullet$ Third step}. 
Recall that we have
\begin{itemize}
\item[(*)]   $\tau_{L,\varphi} \tau_\varphi + \displaystyle\sum_{X\in K_G(H)} \gamma_{L,X} = 0$
\item[(**)]  $\tau_{\varphi,L}\gamma_{L,X}=\tau_\varphi\gamma_X$
\item[(***)] $\tau_{\varphi,L}\tau_{L,\varphi} + \displaystyle\sum_{X\in K_G(H)} \gamma_X = 0$
\end{itemize}
Left-multiplying (*) by $\tau_{\varphi,L}$ and right-multiplying (***) by $\tau_\varphi$ we get:
$$(\sum_{X\in K_G(H)} \gamma_X) \tau_\varphi =
               \tau_{\varphi,L} (\sum_{X\in K_G(H)} \gamma_{L,X})$$
Then by (**) we replace the right term to get 
$$(\sum_{X\in K_G(H)} \gamma_X) \tau_\varphi =\tau_\varphi(\sum_{X\in K_G(H)}\gamma_X)$$
hence $(\sum_{X\in K_G(H)}\limits \gamma_X)$ and $\tau_\varphi$ commute. If $G$ is cyclic, the sum has only one term $\gamma_X$ commuting with (the unique) $\tau_\varphi$.
\findemo

%%%%%%%%%%%%%%%%%%%%%%%%%%%%%%%%%%%%%%%%%%%%%%%%%%%%%%%
\section{Inductive construction of a basis of~$\Ee$}%%%
%%%%%%%%%%%%%%%%%%%%%%%%%%%%%%%%%%%%%%%%%%%%%%%%%%%%%%%
Recall that $k$ be a field of characteristic $p\ge 3$ and $G$ is an elementary abelian $p$-group of rank $r$.
Our aim is to construct a basis of the algebra $\Ext^*_{\comack_k(G)}(S_\un^G,S_\un^G)$. We do this by induction. The technical structure of the induction is build on the following general result on commutative diagrams~:
\begin{mth}{Lemma}\label{diagrchase}
Suppose we have the following diagram of finite dimensional $k$-vector spaces, where the four exterior triangles and the triangle $(E_0,F_2,F_1)$ are commutative. 
$$
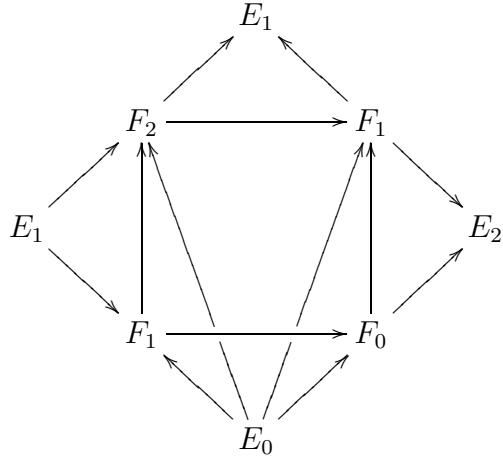
\xymatrix{
  && E_1 &&\\
  & F_2\ar[rr]\ar[ru] && F_1\ar[lu]\ar[rd] & \\
  E_1\ar[ru]\ar[rd] &&&& E_2 \\
  & F_1\ar[uu]\ar[rr] && F_0\ar[uu]\ar[ru] & \\
  && E_0\ar[lu]\ar[ru]\ar[luuu] |!{[lu];[ru]}\hole \ar[ruuu] |!{[lu];[ru]}\hole &&}
$$
Suppose that, in the above diagram, the two maps $E_0\to F_1$ are the same and that all the sequences
$$\xymatrix{
E_r\ar[r]&F_i\ar[r]^{\alpha}&F_j\ar[r]&E_s\,,
}$$
are exact, where $\alpha$ is any of the horizontal or vertical maps. 
Suppose moreover that
\begin{itemize}
\item[(H1)] $\Ker(E_0\to F_0)=\Ker(E_0\to F_1)$,
\item[(H2)] $\Ker(E_1\to F_1)=\Ker(E_1\to F_2)$,
\item[(H3)] $E_0\to F_0\to E_2$ is exact,
\item[(H4)] $\Img(E_1\to F_1)=\Img(F_0\to F_1)$.
\end{itemize}
Then we have 
\begin{itemize}
\item[(C1)] $\Img(F_1\to E_1)=\Img(F_2\to E_1)$
\item[(C2)] $\Img(F_0\to E_2)=\Img(F_1\to E_2)$
\item[(C3)] $E_0\to F_1\to E_2$ is exact
\item[(C4)] $\Img(E_1\to F_2)=\Img(F_1\to F_2)$ and $F_0\cong F_1\cong F_2$
\item[(C5)] $E_1\to F_2\to E_1$ is exact. 
\end{itemize}
\end{mth}

\pf
Choose a direct decomposition $F_0=U\oplus V$ such that 
$$E_0\cong \Ker(E_0\to F_0)\oplus U\;\;.$$ 
Given (H1) and the fact that the triangle $(E_0,F_1,F_0)$ is commutative, the horizontal map in this triangle induces an isomorphism between $\Img(E_0\to F_1)$ and $\Img(E_0\to F_0)$. Similarly, (H2) and the fact that the triangle $(E_1,F_1,F_2)$ is commutative, implies that the vertical map in this triangle induces an isomorphism between $\Img(E_1\to F_1)$ and $\Img(E_1\to F_2)$. \par
Moreover, one has $\Img(E_0\to F_0)\cong U$ and, using (H4) and the exactness of the sequence $E_0\to F_0\to F_1$, that $\Img(E_1\to F_1)\cong V$.

Moreover $\Img(F_1\to F_0)=\Ker(F_0\to E_2)=\Img(E_0\to F_0)\cong U$. 
Thus we have $F_1\cong\Img(E_1\to F_1)\oplus \Img(F_1\to F_0)\cong V\oplus U$, and it follows that $\Img(F_1\to F_2)\cong F_1/\Img(E_0\to F_1)\cong V$ and 
$$\Img(F_1\to E_1)\cong F_1/\Img(F_0\to F_1)\cong (U\oplus V)/V\cong U\;\;.$$
Since the sequence $F_1\to F_2 \to E_1$ is exact, one gets 
$$F_2\cong V\oplus \Img(F_2\to E_1)\;\;.$$ 
By the commutativity of the triangle $(E_1,F_2,F_1)$,
 $$\Img(F_2\to E_1)\subseteq\Img(F_1\to E_1)\cong U\;\;.$$ 
Since the sequence $E_1\to F_2\to F_1$ is exact,
we get that 
$$F_2\cong V\oplus\Img(F_2\to F_1)\;\;.$$ 
By the commutativity of the triangle $(E_0,F_2,F_1)$, we have 
$$\Img(F_2\to F_1)\supseteq \Img(E_0\to F_1)\cong U\;\;.$$ 
Hence $F_2$ is contained in and, respectively contains a $k$-vector space isomorphic to $U\oplus V$ implying that $F_2\cong U\oplus V$ and that the previous inclusions are equalities. \par
In particular $\Img(F_2\to E_1)\cong U$, hence (C1). Moreover, $\Img(E_1\to F_2)=\Img(F_1\to F_2)\cong V$, and (C4) follows. 
Also 
$$\Img(F_1\to E_2)\cong F_1/\Img(F_2\to F_1)\cong (U\oplus V)/U\cong V\;\;,$$ 
and, using (H3), we have $\Img(F_0\to E_2)\cong F_0/\Img(E_0\to F_0)\cong V$. The commutativity of the triangle $(E_2,F_0,F_1)$ implies then (C2). 
In the triangle $(E_0,F_1,F_2)$, we have $\Img(E_0\to F_1)=\Img(F_2\to F_1)$ and, from the exactness of the sequence $F_2\to F_1\to E_2$, we have $\Img(F_2\to F_1)=\Ker(F_1\to E_2)$. Now (C3) follows. Lastly, $\Img(E_1\to F_2)=\Img(F_1\to F_2)=\Ker(F_2\to E_1)$ and (C5) follows. 
\findemo

\npar\label{LongExactSequences} First we fix some notation~: we set $F_a(n):=\Ext^n_{\comack_k(G)}(T_a,S_\un^G)$. In particular, $E_G(n)=F_1(n)$ and $L(n)=F_{p-2}(n)$. Applying $\Hom_{\comack_k(G)}(-,S_\un^G)$ to the short exact sequences of Mackey functors in~\ref{kTaTa+1}, yields two long exact sequences
$$0\to E_G(0)\to F_{a+1}(0)\to F_a(0)\to E_G(1)\to \dots$$
and
$$0\to F_{a}(0)\to F_{a+1}(0)\to E_G(0)\to F_{a}(1)\to \dots$$
In what follows $E_G(n-1)\to F_a(n)$, $E_G(n)\to F_a(n)$, $F_a(n)\to E_G(n)$, $F_a(n)\to E_G(n+1)$, $F_a(n)\to F_{a+1}(n)$ or $F_{a+1}(n)\to F_{a}(n)$ are the maps in the long exact sequences above.
\begin{mth}{Proposition}\label{prop:L(n)=E(n)}
Let $p>5$. Then, for all $n$, the sequence 
$$E_G(n)\to L(n)\to E_G(n)$$ 
given as above is exact and $\dim L(n)=\dim E_G(n)$. Moreover, the sequence 
$$E_G(n)\to E_G(n+1)\to E_G(n+2)$$ 
given by multiplication from the right with $\tau_\varphi$, is exact, for any non trivial homomorphism $\varphi:G\to k^+$ and we have a direct decomposition of $E_G(n)$ inductively given by
$$\begin{array}{lllll} \vspace{3pt}
E_G(n) &\cong\,\, &E_H^G(n) &\oplus \,\,&E_H^G(n-1)\tau_\varphi \\ \vspace{3pt}
\qquad &\oplus &E_H^G(n-2)s_\varphi &\oplus &E_H^G(n-3)\tau_\varphi s_\varphi \\ \vspace{3pt}
\qquad &\oplus &E_H^G(n-4)s_\varphi^2 &\oplus &E_H^G(n-5)\tau_\varphi s_\varphi^2 \\ \vspace{3pt}
\qquad &\oplus &\dots &&\\ \vspace{3pt}
\qquad &\oplus &\displaystyle\bigoplus_{X\in\mathcal{X}}E_G(n-2) \gamma_X &\oplus &\displaystyle\bigoplus_{X\in\mathcal{X}}E_G(n-3) \gamma_X\tau_\varphi \\ \vspace{2pt}
\qquad &\oplus &\displaystyle\bigoplus_{X\in\mathcal{X}}E_G(n-4) s_\varphi\gamma_X &\oplus &\displaystyle\bigoplus_{X\in\mathcal{X}}E_G(n-5) s_\varphi\gamma_X \tau_\varphi \\ \vspace{2pt}
\qquad &\oplus &\displaystyle\bigoplus_{X\in\mathcal{X}}E_G(n-6) s_\varphi^2\gamma_X &\oplus &\displaystyle\bigoplus_{X\in\mathcal{X}}E_G(n-7) s_\varphi^2\gamma_X \tau_\varphi \\ \vspace{2pt}
\qquad &\oplus &\dots &&\\
\end{array}$$
where $E_H^G(m):=\sigma_{G/Y}^G\Iso_{H}^{G/Y} E_H(m)$ and $s_\varphi:=\sum_{X\in K_G(H)}\limits \gamma_X$.
\end{mth}

\npar\label{StartInduction} We do the proof by simultaneous induction on $n$ and $|G|$. The cases $n=0$ or $G=\un$ are trivial so we have the starting point for the induction. We suppose the proposition true for all the proper subgroups $K$ of $G$ and all $m$, and, respectively, for $K=G$ and $m\le n$, and we prove it for $G$ and $n+1$.  

\npar By induction, we have that the sequence 
$$E_K(m-1)\to E_K(m)\to E_K(m+1)$$ 
given by multiplication from the right with $\tau_\psi$, is exact, for $K$ and $m$ as in the previous paragraph, and for any non trivial homomorphism $\psi:K\to k^+$. To ease notation, we set $E(n):=E_G(n)$, while still writing $E_K(n)$ when $K$ is different from $G$.

%Also by induction, we have that $\Ker\big(E(n-1)\to F_a(n)\big)=\Ker\big(E(n-1)\to F_b(n)\big)$ for all $1\le a,b\le p-3$ and that $\Ker\big(E(n)\to F_a(n)\big)=\Ker\big(E(n)\to F_b(n)\big)$ for all $2\le a,b\le p-2$.

\npar The four exterior triangles and the triangle $\big(E(n-1),F_{d+2}(n),F_{d+1}(n)\big)$ in the following diagram
$$
\xymatrix{
  && E(n) &&\\
  & F_{d+2}(n)\ar[rr]\ar[ru] && F_{d+1}(n)\ar[lu]\ar[rd] & \\
  E(n)\ar[ru]\ar[rd] &&&& E(n+1) \\
  & F_{d+1}(n)\ar[uu]\ar[rr] && F_{d}(n)\ar[uu]\ar[ru] & \\
  && E(n-1)\ar[lu]\ar[ru]\ar[luuu] |!{[lu];[ru]}\hole \ar[ruuu] |!{[lu];[ru]}\hole  &&}
$$
are commutative for all $n$ and $1\le d\le p-4$. Setting $E_i=E(n+i-1)$ and $F_i=F_{d+i}(n)$ for $i\in\{0,1,2\}$, for $d=1$ the corresponding hypothesis (H3) in Lemma~\ref{diagrchase} is the fact that $E(n-1)\to E(n)\to E(n+1)$ is exact and the hypotheses (H1), (H2), (H4) are easy to check. This starts the induction on~$d$. Then, (H1) to (H4) are satisfied for all $d\le p-4$. Indeed, the conclusions (C3) and (C4) for $F_i=F_{d+i}(n)$ in the same Lemma~\ref{diagrchase} are the hypotheses (H3) and (H4) for $F_i=F_{d+i+1}(n)$. Also, we have
$$\Ker\big(E(n)\to F_a(n+1)\big)=\Ker\big(E(n)\to F_b(n+1)\big),\, \forall a,b\in\{ 1,\ldots, p-3\}\text{ and}$$ 
$$\Ker\big(E(n+1)\to F_a(n+1)\big)=\Ker\big(E(n+1)\to F_b(n+1)\big),\, \forall a,b\in\{ 2,\ldots,p-2\}.$$
Hence (H1)-(H2) for $F_i=F_{i+d+1}(n)$ are obtained from (C1) and (C2) for $F_i=F_{d+i}(n)$ and the long exact sequences in~\ref{LongExactSequences}: 
$$\Img\big(F_a(n)\to E(n)\big)=\Ker\big(E(n)\to F_{a-1}(n+1)\big),\, \forall a\in\{2,\ldots, p-2\}\text{ and}$$ 
$$\Img\big(F_a(n)\to E(n+1)\big)=\Ker\big(E(n+1)\to F_{a+1}(n+1)\big), \forall a\in\{1,\ldots,p-3\}\,.$$

\npar Thus, by induction on~$d$ and using Lemma~\ref{diagrchase}, (C4) and (C5) are true for all $d\in\{1,\ldots,p-4\}$. In particular, (C5) for $d=p-4$, in which case $F_2=F_{p-2}(n)=L(n)$, together with the last part of (C4), for any integer $d\in\{1,\ldots,p-4\}$, give the first two claims in Proposition~\ref{prop:L(n)=E(n)}.

\npar The last ingredient we need to continue the induction on $n$ is that the sequence $E(n)\to E(n+1)\to E(n+2)$ is exact. These maps are given by multiplication from the right by $\tau_\varphi$ in Yoneda's notation. To prove the exactness of the sequence at $E(n+1)$, we fix a decomposition $G=H\oplus Y$, where $Y$ has order $p$, and we explicitly decompose $E(n+1)$ with respect to $E(m)$ for $m\le n$ and to $E_H^G(m):=\sigma_{G/Y}^G\Iso_{H}^{G/Y} E_H(m)$ for $m\le n+1$.

\npar Recall that we have the following properties:
\begin{itemize}
\item $E(n+1) \cong
L(n)\tau_{L,\varphi} \oplus \displaystyle\bigoplus_{X\in\mathcal{X}}E(n-1)\gamma_X \oplus E_H^G(n+1)$, where we have set $\mathcal{X}=K_G(H)\setminus\{Y\}$.
\item the multiplication by $\gamma_X$ induces an injective map from $E(m)$ to ${E(m+2)}$.
\item $E(n)\to L(n)\to E(n)$ is exact and $\dim L(n)=\dim E(n)$.
\item the multiplication by $\tau_{L,\varphi}$ induces an injective map from $L(m-1)$ to $E(m)$.
\end{itemize}

\npar Let's concentrate on $L(n)\tau_{L,\varphi}$. It is easy to check that 
$$\Img\big(E(n)\to L(n)\big)\tau_{L,\varphi}=E(n)\tau_\varphi\;\;.$$ 
Moreover from the last diagram of the induction on $d$ we have 
$$L(n)=\Img\big(E(n-1)\to L(n)\big)+\Img\big(E(n)\to L(n)\big)\;\;.$$ 
In terms of Yoneda's composition of extensions, this becomes 
$$ L(n)\tau_{L,\varphi} = E(n)\tau_\varphi + E(n-1)\tau_{\varphi,L}\tau_{L,\varphi}=
   E(n)\tau_\varphi + E(n-1)(\sum_{X\in K_G(H)} \gamma_X)\,.$$
The above sum is not direct but we have a description of the intersection of the two terms. This is given in the following proposition.
\pagebreak[3]
\begin{mth}{Proposition}
\begin{eqnarray*}
(\#)\quad &:=&E(n)\tau_\varphi \cap E(n-1)s_\varphi \gamma_X \\
\qquad    &= &E(n-2) \tau_\varphi s_\varphi \\
\qquad    &=& E(n-2) s_\varphi \tau_\varphi
\end{eqnarray*}
\end{mth}
\pf
By the induction started in~\ref{StartInduction}, the kernel of the multiplication from the left by $\tau_\varphi$ from $E(n)$ to $E(n+1)$ is $E(n-1)\tau_\varphi$. Hence $\dim (E(n-2)s_\varphi\tau_\varphi) =
\dim\big(E(n-2)s_\varphi\big) - \dim \big(E(n-1)\tau_\varphi \cap E(n-2)s_\varphi\big)$

Also by induction on $n$ one can assume that 
$E(n-1)\tau_\varphi \cap E(n-2)s_\varphi = E(n-3)\tau_\varphi s_\varphi$ and that the 
multiplication by $s_\varphi$ gives an injective map. Thus 
$$\begin{array}{ll}
\dim\big(E(n-2)s_\varphi\tau_\varphi\big)
            &= \dim \big(E(n-2)s_\varphi\big) - \dim \big(E(n-3)\tau_\varphi s_\varphi\big)\\
            &= \dim E(n-2) - \dim \big(E(n-3)\tau_\varphi\big)\\
            &= \dim E(n-2)\tau_\varphi
\end{array}$$          
where the last equality is given by the exactness of the sequence
$$E(n-3)\to E(n-2)\to E(n-1)$$ 
with maps given by multiplication from the right by $\tau_\varphi$.
Thus we get that the multiplication from the right by $s_\varphi$ gives a bijection
between $E(n-2)\tau_\varphi$ and $E(n-2)\tau_\varphi s_\varphi$.

Moreover $E(n-2)\tau_\varphi s_\varphi$ is contained in the intersection $(\#)$ which enables the following computation:
\begin{eqnarray*}
\dim \big(E(n-2)\tau_\varphi s_\varphi\big) &\le& \dim(\#)\\
            &=&\dim E(n)\tau_\varphi + \dim E(n-1)s_\varphi - \dim L(n)\tau_{L,\varphi} \\
            &=&\dim E(n)\tau_\varphi + \dim E(n-1)s_\varphi - \dim L(n) \\
            &=&\dim E(n-1)s_\varphi - \dim E(n-1)\tau_\varphi \\
            &\le&\dim E(n-1) - \dim E(n-1)\tau_\varphi \\
            &=&\dim E(n-2)\tau_\varphi
\end{eqnarray*}         
We proved that the first and the last expression in the sequence are equal so all the inequalities are equalities and we get $(\#)$ and, moreover, that $\dim E(n-1) = \dim E(n-1)s_\varphi$.
\findemo

\npar Making use of the previous two propositions we get the following decomposition of $E(n+1)$ as a direct sum~:
$$
E(n+1)\cong [E(n)\tau_\varphi + E(n-1)(\sum_{X\in K_G(H)} \gamma_X)] \oplus\displaystyle\bigoplus_{X\in\mathcal{X}}E(n-1)\gamma_X \oplus E_H^G(n+1).
$$
\npar\label{decomposition} Using this, by induction on $n$ one can get the following direct sum decomposition, which has the advantage that for every term, the multiplications from the right by monomials of type $s_\varphi^\alpha\gamma_X^\delta\tau_\varphi^\varepsilon$, for $\alpha\ge 0$, for $\delta,\varepsilon\in\{0,1\}$, and $(\alpha,\delta)\ne(0,0)$, induce bijective maps from $E(n+1-2\alpha-2\delta-\varepsilon)$ to $E(n+1)$~:
$$\begin{array}{lllll} \vspace{3pt}
E(n+1) &\cong\,\, &E_H^G(n+1) &\oplus \,\,&E_H^G(n)\tau_\varphi \\ \vspace{3pt}
\qquad &\oplus &E_H^G(n-1)s_\varphi &\oplus &E_H^G(n-2)\tau_\varphi s_\varphi \\ \vspace{3pt}
\qquad &\oplus &E_H^G(n-3)s_\varphi^2 &\oplus &E_H^G(n-4)\tau_\varphi s_\varphi^2 \\ \vspace{3pt}
\qquad &\oplus &\dots &&\\ \vspace{3pt}
\qquad &\oplus &\displaystyle\bigoplus_{X\in\mathcal{X}}E(n-1) \gamma_X &\oplus &\displaystyle\bigoplus_{X\in\mathcal{X}}E(n-2) \gamma_X\tau_\varphi \\ \vspace{2pt}
\qquad &\oplus &\displaystyle\bigoplus_{X\in\mathcal{X}}E(n-3) s_\varphi\gamma_X &\oplus &\displaystyle\bigoplus_{X\in\mathcal{X}}E(n-4) s_\varphi\gamma_X \tau_\varphi \\ \vspace{2pt}
\qquad &\oplus &\displaystyle\bigoplus_{X\in\mathcal{X}}E(n-5) s_\varphi^2\gamma_X &\oplus &\displaystyle\bigoplus_{X\in\mathcal{X}}E(n-6) s_\varphi^2\gamma_X \tau_\varphi \\ \vspace{2pt}
\qquad &\oplus &\dots &&\\
\end{array}$$

\npar Now, for every group homomorphism $\psi:G\to k^+$ with $\Ker\,\psi\ne H=\Ker\,\varphi$, the restriction from $E(m)$ to $E_H(m)$ gives an isomorphism between the multiplication by $\tau_\psi$ on $E_H^G(m)$ and the multiplication by $\tau_{\Res_H^G(\psi)}$.  Hence  the short exact sequence $E_H(m)\to E_H(m+1)\to E_H(m+2)$ given by multiplication from the left by $\tau_{\Res_H^G(\psi)}$ induces a short exact sequence 
$$E_H^G(m)\to E_H^G(m+1)\to E_H^G(m+2)$$ 
given by the multiplication from the left by $\tau_\psi$. Summing all up we get an exact sequence $E(n)\to E(n+1)\to E(n+2)$ given by the multiplication from the left by~$\tau_\psi$.

\npar The duality sends the multiplication from the left by $\tau_\psi$ to multiplication from the right by $\tau_\psi$ and we can exchange the roles of $\varphi$ and $\psi$ to get that the exact sequence $E(n)\to E(n+1)\to E(n+2)$ given by multiplication from the right by $\tau_\varphi$ is exact for any homomorphism $\varphi:G\to k^+$. This finishes the induction step on $n$ and the proof of Proposition~\ref{prop:L(n)=E(n)}.

\npar\label{decomposition for p=3} Remark that for $p=3$ we have $\tau_\varphi^2=s_\varphi$ thus, also in this case, $\tau_\varphi$ and $s_\varphi$ commute. Now, using in Corollary~\ref{short exact sequence} the fact that $L(n)=E(n)$ we get the following decomposition as a direct sum 
$$E(n+1)\cong E(n)\tau_\varphi\oplus\bigoplus_{X\in\mathcal{X}}E(n-1)\gamma_X\oplus E_H^G(n+1)\mpoint$$
Decomposing $E(n)$ in the same way we get
$$E(n+1)\!\cong\! E(n-1)s_\varphi\oplus\bigoplus_{X\in\mathcal{X}}\!E(n-2)\gamma_X\tau_\varphi \oplus E_H^G(n)\tau_\varphi\oplus\bigoplus_{X\in\mathcal{X}}\!E(n-1)\gamma_X\oplus E_H^G(n+1)\,$$
Continuing the replacement, we get the same decomposition as in~\ref{decomposition}.

\npar Following the proof of Theorem 14.2 in \cite{cohocplx} we get now the Poincar\'e series for $\Ext^*_{\comack_k(G)}(S_\un^G,S_\un^G)$ for all odd $p$. Indeed, the only ingredients needed in that proof are the exactness of 
$$0\longrightarrow L(n-1)\oplus\mathop{\oplus}_{X\in\mathcal{X}}\limits E_G(n-2)\to E_G(n)\to E_H(n)\to 0\mpoint$$
and that $\dim E_G(n)=\dim L(n)$. These facts are proved in Corollary~\ref{short exact sequence}, respectively in Proposition~\ref{prop:L(n)=E(n)}.

%%%%%%%%%%%%%%%%%%%%%%%%%%%%%%%%%%%%%%%%%%%%%%
\section{A presentation of the algebra $\Ee$}%%%
%%%%%%%%%%%%%%%%%%%%%%%%%%%%%%%%%%%%%%%%%%%%%%
We are now able to give a presentation of the algebra $\Ee=\Ext^*_{\comack_k(G)}(S_\un^G,S_\un^G)$. Our aim is to present an algbra $\A$ that, a priori, has $\Ee$ as quotient and then to show that the two algebra are isomorphic. Let $r$ be the rank of $G$ and ${0}=H_0<H_1<\cdots<H_{r-1}=H<H_r=G$ be a maximal flag in $G$. We choose a direct sum decomposition of $G$
$$G=Y_1\oplus Y_2\oplus\cdots\oplus Y_{r-1}\oplus Y_r\mvirg$$ 
where $Y_i$ is a complement of $H_{i-1}$ in $H_i$. Let $\varphi_i:G\to \F_p$ be a non-trivial morphism with kernel $\oplusb{1\le j\le r}{j\ne i}\limits Y_j$. For every $i\in\{1,\ldots,r\}$ we define an atom $\hat\tau_i$ of degree $1$, and for every $i\in\{1,\ldots,r\}$ and for every subgroup $X\le G$ of order $p$,  such that $XH_{i-1}=H_i$, we say that $X$ {\em has position i}, and we define an atom $\hat\gamma_X$ of degree $2$. \par
Consider the set $A=\{\hat\tau_i|1\le i\le r\}\cup\{\hat\gamma_X|X\le G,\, |X|=p\}$. A {\it word} with atoms in~$A$ is obtained in the natural way, by concatenation. The {\it degree} of a word is the sum of the degrees of its atoms. Denote by $\rk X$ the position of~$X$. We say that a word is {\it of type} $i$ if the atoms appearing in its decomposition are among $\hat\tau_i$ and $\hat\gamma_X$ with $\rk X=i$. The empty word is of type $i$ for all $i\in\{1,\ldots,r\}$.
\npar For all $n\ge 0$, let $A_n$ be the set of words of degree $n$ with the property that they don't contain any of the following sequences of two consecutive atoms:
\begin{itemize}
\item [(S1)]$\hat\tau_i\hat\tau_j$, for $1\le j\le i\le r$;
\item [(S2)]$\hat\tau_i\hat\gamma_X$, for $1\le\rk X<i\le r$;
\item [(S3)]$\hat\gamma_{Y_i}\hat\tau_i$, for $1\le i\le r$;
\item [(S4)]$\hat\gamma_X\hat\gamma_Y$, for $1\le \rk Y< \rk X\le r$;
\item [(S5)]$\hat\gamma_X\hat\tau_i$, for $1\le i <\rk X\le r$.
\end{itemize}
We call {\it pre-admissible} the words $w$ that can be written as a product $w=w_1w_2\dots w_r$ where, for $1\le i\le r$, the word $w_i$ is of type $i$. We call {\it admissible} the pre-admissible words that do not contain any of the sequences $\hat\tau_i\hat\tau_i$ and $\hat\gamma_{Y_i}\hat\tau_i$.

\begin{mth}{Proposition}\label{CanonicalForm}
Let $w$ be a word of degree $n$. Then $w$ belongs to the union of the $A_n$'s for all $n$ if and only if it is admissible.
\end{mth}
\pf
The forbidden sequences (S1) for $i\ne j$, (S2), (S4) and (S5) imply that for $j<i$, the atoms $\hat\tau_i$ and $\hat\gamma_X$ with $\rk X=i$ only appear after $\hat\tau_j$ and $\hat\gamma_Y$ with $\rk Y=j$. This gives a pre-admissible form $w_1w_2\dots w_r$ for $w$. Moreover, in each $w_i$, the sequences $\hat\tau_i\hat\tau_i$ and $\hat\gamma_{Y_i}\hat\tau_i$ are explicitly forbidden by (S1) for $i=j$ and (S3). 
\findemo

\begin{mth}{Proposition}\label{DecompBasis}
Let $p$ be an odd prime and $\A$ be the graded $k$-algebra with generators 
$$\{\hat\tau_i|1\le i\le r\}\;\hbox{in degree 1, and }\{\hat\gamma_X|X\le G,\, |X|=p\}\;\hbox{in degree 2},$$ 
subject to the relations 
\begin{itemize}
\item [(R1)]$\hat\tau_i\hat\tau_i=0$, if $p\geq 5$, or \\
$\hat\tau_i\hat\tau_i=-\sum_{X\nleq\Ker\,\varphi_i}\limits\hat\gamma_X$ if $p=3$.
\item [(R2)]$\hat\tau_i\hat\tau_j+\hat\tau_j\hat\tau_i=0$ for $1\le j\le i \le r$, if $p\geq 5$, or \\
$\hat\tau_i\hat\tau_j+\hat\tau_j\hat\tau_i=\sum_{X\nleq\Ker(\varphi_i+\varphi_j)}\limits\hat\gamma_X-\sum_{X\nleq\Ker\,\varphi_i}\limits\hat\gamma_X-\sum_{X\nleq\Ker\,\varphi_j}\limits\hat\gamma_X$ if $p=3$.
\item [(R3)]$[\hat\tau_i,\sum_{X\nleq\Ker\,\varphi_i}\limits\hat\gamma_X]=0$;
\item [(R4)]$[\varphi_j(x)\hat\tau_i-\varphi_i(x)\hat\tau_j,\hat\gamma_{\langle x\rangle}]=0, \text{ for }1\le i<j\le 1,\,x\in G$;
\item [(R5)]$[\hat\gamma_X,\sum_{Y<Q}\limits\hat\gamma_Y]=0,\text{ for all }\,X<Q\le G,\,|X|=p,\,|Q|=p^2$.
\end{itemize}
Then for all $n\ge 1$, $A_n$ is a basis of the $k$-vector space of elements of degree~$n$ in~$\A$.
\end{mth}
\begin{rem}{Remark}
In the case $p=3$, Relation~(R3) is implied by Relation~(R1). The relations (R1) to (R5) are rather technical. This comes from the need to have a finite, easy to order set of generators, appropriate for induction arguments. We can give a more intuitive presentation for the algebra $\A$ if we allow ourselves to increase the set of generators, by taking a generator $\hat\tau_\varphi$ for every group homomorphism $\varphi:G\to \F_p^+$. Here we identify $\hat\tau_{\varphi_i}$ with $\hat\tau_i$. The relations are then given by
\begin{breakbox}
\begin{itemize}
\item [(L1)]$\hat\tau_\varphi+\hat\tau_\psi=\hat\tau_{\varphi+\psi}$.
\item [(L2)]$\hat\tau_\varphi\hat\tau_\varphi=0$ and $[\hat\tau_\varphi,\sum_{X\nleq\Ker\,\varphi}\limits\hat\gamma_X]=0$, if $p\geq 5$, or \\
$\hat\tau_\varphi\hat\tau_\varphi=-\sum_{X\nleq\Ker\,\varphi}\limits\hat\gamma_X$ if $p=3$.
\item [(L3)]$[\hat\tau_\varphi,\hat\gamma_X]=0, \text{ for }X\le\Ker\,\varphi,|X|=p$.
\item [(L4)]$[\hat\gamma_X,\sum_{Y<Q}\limits\hat\gamma_Y]=0,\text{ for all }\,X<Q\le G,\,|X|=p,\,|Q|=p^2$.
\end{itemize}
\end{breakbox}
We leave it as an easy exercise to the reader to verify that the two presentation lead to isomorphic algebras. These two presentations yield Theorem~\ref{presentation Ext algebra} and~\ref{presentation Ext algebra 1}, respectively.
\end{rem}
To prove Proposition~\ref{DecompBasis}, we need the following sequence of lemmas.

\begin{mth}{Lemma}\label{CommuteGammaXGammaY}
Let $l<m$ and let $X$ and $Y$ be subgroups of order $p$ of $G$, such that $\rk X=m$ and $\rk Y=l$. Then we have $\hat\gamma_X\hat\gamma_Y= \hat\gamma_Y\hat\gamma_X + W$ where $W$ is a sum of admissible words $\hat\gamma_{X'}\hat\gamma_{X''}$ of type $m$ and degree~$4$.
\end{mth}
\pf
Firstly we prove that all the subgroups of order $p$ of $\langle X,Y\rangle$, different from $Y$, have position $m$. Indeed, if $x$ is a generator of $X$ and $y$ is a generator of $Y$, then denote by $X_c$ the subgroup of $\langle X,Y\rangle$ generated by $x+cy$ for $1\le c \le p-1$. As $cy\in Y\le H_{m-1}$ and $XH_{m-1}=H_m$ we get that $X_cH_{m-1}=H_m$ so $X_c$ has position $m$. 

Secondly we use Relation (R5) with $Q=\langle X,Y\rangle$ to get $$\hat\gamma_X\hat\gamma_Y+\displaystyle\sum_{c=1}^{p-1}\hat\gamma_X\hat\gamma_{X_c}=\displaystyle\sum_{c=1}^{p-1}\hat\gamma_{X_c}\hat\gamma_X+\hat\gamma_Y\hat\gamma_X$$
Besides $\hat\gamma_X\hat\gamma_Y$, all the words appearing in the above equation are admissible and, excepting $\hat\gamma_Y\hat\gamma_X$, they are all of type $m$. The result follows. 
\findemo

\begin{mth}{Lemma}\label{CommuteGammaXTaul}
Let $l<m$ and $X$ of order $p$ with $\rk X=m$. Then $\hat\gamma_{X}\hat\tau_l=\hat\tau_l\hat\gamma_X+W$, where $W$ is a linear combination of pre-admissible words of type $m$ and degree $3$.
\end{mth}
\pf
Let $x$ be a generator of $X$. If we set $\alpha:=\varphi_l(x)$ and $\beta:=\varphi_m(x)$, relation (R4) gives that $$\hat\gamma_X(\beta\hat\tau_l-\alpha\hat\tau_m)=(\beta\hat\tau_l-\alpha\hat\tau_m)\hat\gamma_X$$
As $\hat\tau_l\hat\gamma_X$, $\hat\gamma_X\hat\tau_m$ and $\hat\tau_m\hat\gamma_X$ are all pre-admissible words of degree $3$, the latter two being of type $m$, and $\beta\ne 0$ the result follows. The fact that we only get pre-admissibility (and not admissibility) comes from the fact that we might have $X=Y_m$.
\findemo

\begin{mth}{Lemma}\label{PreCanonicalForm}
Every word in $w\in\A$ can be written as a finite sum of pre-admissible words, i.e. $w=\sum_{j=1}^s w_{j1}w_{j2}\dots w_{jr}$, with $w_{ji}$ of type $i$ for $1\le i\le r$ (the $w_{ji}$'s are allowed to be the empty word). Moreover, if $w$ is composed of atoms of position at least $l$, then the $w_{ji}$'s are empty for all $i<l$ and $j=1,2,\dots,s$.  
\end{mth}
\pf
We define a quadruple lexicographic order $\Oo$ on the set of words in $\A$ first by their degree then, for equal degrees, by the number of atoms of degree $1$ (i.e. the $\hat\tau_i$'s) appearing in the word, then by the degree of the left term in the decomposition $uv$, where $v$ is pre-admissible, of maximal possible degree, then, finally, by the number of atoms in $v$ of higher position than the rightmost atom in $u$. We prove the lemma by induction on the order $\Oo$. We'll see in the proof that the order by the number of atoms of degree $1$ appearing in the word is only needed for the case $p=3$. 

The lemma is trivial for words of degree $1$ and straight forward, using Relations (R1) and (R2), for words of degree $2$. Now take a word $w\in\A$ of degree $d\ge 3$. Write $w=uv$, a concatenation, with $v$ pre-admissible, of maximal possible degree. If $\deg u=0$, i.e. $u$ is empty, then $w$ is pre-admissible and we are done. If the number of atoms in $v$ of higher position than the rightmost atom in $u$ is $0$ or the position of the leftmost atom in $v$ is $n$ then $w=v$ and we are again done. It is very important to remark that all the manipulations we make in the proof, using the relations in Proposition~\ref{DecompBasis}, are not increasing the number of atoms of order $1$. Suppose now the lemma is true for all words smaller than $w$ for the order $\Oo$. Let $a$ be the rightmost atom of $u$ i.e. $u=w'a$. If $u\ne a$ then $av$ is of degree smaller than $w$ so by the induction hypothesis on the degree $av= \sum_{j=1}^s w_j$ with $w_j$ pre-admissible. Moreover $w'w_j$ is s
 maller than $w$ with respect to the order $\Oo$, for all $j=1,2,\dots,s$, so, by the induction hypothesis, $w'w_j$ is a sum of pre-admissible  elements. Also by the induction hypothesis, the atoms appearing in the new words are of higher or equal position than the ones appearing in $w$. Thus we can suppose that $u=a$ is an atom of position $m$ so $u=\hat\tau_m$ or $u=\hat\gamma_X$ with $\rk X=m$.

Case 1. Let $u=\hat\gamma_X$ with $X$ of position $m$. We distinguish two sub-cases:
\begin{itemize}
\item[i)] $v=\hat\gamma_Yv'$ with $Y$ of position $l$ and $l<m$. Then, using Lemma~\ref{CommuteGammaXGammaY}, we have $\hat\gamma_Xv=\hat\gamma_X\hat\gamma_Yv'=\hat\gamma_Y\hat\gamma_Xv'+Wv'$. We have that $\hat\gamma_Xv'$ is smaller than $w$, hence, by the induction hypothesis, it is equal to a sum of pre-admissible words, composed of atoms of position at most $l$. Thus the words obtained from these by multiplying to the left by $\hat\gamma_Y$ are still pre-admissible. The expression $Wv'$ is a sum of words of type $\hat\gamma_{X'}\hat\gamma_{X''}v'$ with the first two atoms of position $m$. Again, by the induction hypothesis $\hat\gamma_{X''}v'$ is equal to a sum of pre-admissible words $w_i$, which, moreover, have less atoms of lower position than $m$ than $\hat\gamma_Yv'$. Thus $\hat\gamma_{X'}w_i$ is smaller than $w$ for the order $\Oo$ and, by the induction hypothesis, the former is equal to a sum of pre-admissible words.  
\item[ii)] $v=\hat\tau_lv'$ with $l<m$. Then, using Lemma~\ref{CommuteGammaXTaul}, there exists a linear combination $W$ of words of type $m$ such that $\hat\gamma_Xv=\hat\gamma_X\hat\tau_lv'=\alpha\hat\tau_l\hat\gamma_Xv'+Wv'$, with $\alpha\in k$. Proceeding analogously to the part i) above we obtain that  both $\hat\tau_l\hat\gamma_Xv'$ and $Wv'$ are equal to linear combinations of pre-admissible words. 
\end{itemize}

Case 2. Let $u=\hat\tau_m$. We distinguish two sub-cases.
\begin{itemize}
\item[i)] $v=\hat\gamma_Yv'$ with $\rk Y=l$ and $l<m$. Then $Y\le\Ker\,\varphi_m$ and, using Relation (R4), we have $\hat\tau_m\hat\gamma_Yv'=\hat\gamma_Y\hat\tau_mv'$. Moreover $\hat\tau_mv'$ is of degree smaller than $w$ so, by the induction hypothesis, $\hat\tau_mv'$ is a linear combination of pre-admissible words having only atoms of position greater or equal than $l$. So multiplying by $\hat\gamma_Y$ to the left still keeps these words pre-admissible.
\item[ii)] $v=\hat\tau_lv'$ with $l<m$. Then, for $p\ge 5$, using Relation (R2) we have $\hat\tau_mv=\hat\tau_m\hat\tau_lv'=-\hat\tau_l\hat\tau_mv'$ and the result as in part i) of Case 2. When $p=3$ the proof is more difficult. Using also Relation (R2), we have  $\hat\tau_mv=\hat\tau_m\hat\tau_lv'=\sum_{X\nleq\Ker(\varphi_m+\varphi_l)}\limits\hat\gamma_Xv'-\sum_{X\nleq\Ker\,\varphi_m}\limits\hat\gamma_Xv'-\sum_{X\nleq\Ker\,\varphi_l}\limits\hat\gamma_Xv'-\hat\tau_l\hat\tau_mv'$. Moreover $\rk X>l$ for all $\hat\gamma_X$ appearing in the right hand side term of this equality. Now $\tau_mv'$ has smaller degree than $w$, so, by the induction hypothesis, it is equal to a linear combination of pre-admissible words. The other words appearing in the right hand side term have all a smaller number number of atoms of degree $1$ than $w$, so, by the induction hypothesis, they are all equal to linear combinations of pre-admissible words. Remark here that Case 1 and Case 2~i) are complet
 ely solving the case of words having at most one atom of degree $1$, without needing the Case 2~ii). So there is no problem in starting the induction process.
\end{itemize}
\findemo

\bigskip
\pf[of Proposition~\ref{DecompBasis}]
Firstly we prove that $A_n$ is a set of generators in degree $n$. Using Lemma~\ref{PreCanonicalForm} we have that any word of degree $n$ can be written as a linear combination of pre-admissible words $w_1w_2\dots w_r$, i.e. with $w_i$ of type $i$ for $1\le i\le r$. Let $w=w_1w_2\dots w_r$ be one of these words. By Lemma~\ref{CanonicalForm}, it remains to show that, possibly using new linear combinations, we also eliminate the expressions $\tau_i\tau_i$ and $\hat\gamma_{Y_i}\hat\tau_i$ that might appear in the words $w_i$. Let $i$ be lowest index such that $w_i$ has forbidden sequences, with the convention that $i=r+1$ is $w$ contains no forbidden sequence. We define a triple lexicographic order $\Oo$ on the words, given by decreasing order on the index $i$ defined previously, then by the number of atoms of position $i$ and, finally, by the decreasing order of the degree of the its admissible prefix of maximal degree. We proceed by induction on $\Oo$. If the degree of the word
  is at most one or there are no atoms of degree $1$ or the word is already admissible, then there is nothing to prove. Else, denote $w_{_{<i}}:=w_1w_2\dots w_{i-1}$ and by $w_{_{>i}}:=w_{i+1}w_{i+2}\dots w_r$. Let $f$ be the leftmost forbidden sequence in $w_i$. There are two kinds of forbidden sequences to consider.
\begin{itemize}
\item[i)] $f=\hat\tau_i\hat\tau_i$. Write $w=u_1fu_2$. Using Relation (R1), we have $f=0$ when $p\ge 5$ and $f=\sum_{X\nleq\Ker\,\varphi_i}\limits\hat\gamma_X$ when $p=3$. So, the case $p\ge 5$ is trivial ,and in the following we study the case $p=3$. Now for $X\nleq\Ker\,\varphi_i$ we have $\rk X\ge i$. If $\rk X=i$ then the word $w_{_{<i}}u_1\hat\gamma_X$ doesn't have any forbidden sequence, and, by induction hypothesis, $\hat\gamma_Xu_2w_{_{>i}}$ is equal to a linear combination of admissible words. Suppose now that $\rk X>i$. Using Lemma~\ref{PreCanonicalForm}, $\hat\gamma_Xu_2w_{_{>i}}$ is equal to a linear combination of pre-admissible words with atoms of position at least $i$. Let $v_i\dots v_r$ be one of these words. Then $w_{_{<i}}u_1v_i\dots v_r$ is a pre-admissible word with less atoms position $i$ than $w$. By the induction hypothesis, $w_{_{<i}}u_1v_i\dots v_r$ is equal to a linear combination of admissible words.
\item[ii)] $f=\hat\gamma_{Y_i}\hat\tau_i$. Write $w=u_1fu_2$. Relation (R3) gives 
$$\hat\gamma_{Y_i}\hat\tau_i=\hat\tau_i\sum_{X\nleq\Ker\,\varphi_i}\limits\hat\gamma_X+\sum_{X\nleq\Ker\,\varphi_i-\{Y_i\}}\limits\hat\gamma_X\hat\tau_i\,.$$
As $H_{i-1}\le\Ker\,\varphi_i$ we have that any $X\nleq\Ker\,\varphi_i$ is of position at least $i$. The words $w_{_{<i}}u_1\hat\tau_i\hat\gamma_{Y_i}$, $w_{_{<i}}u_1\hat\tau_i\hat\gamma_X$ and $w_{_{<i}}u_1\hat\gamma_X\hat\tau_i$ with $X$ of position $i$ contain no forbidden sequence and are of degree bigger than $w_{_{<i}}u_1$, so, by the induction hypothesis these words multiplied by $w_{_{>i}}$ are equal to linear combinations of admissible words. In the case of the $u_1\hat\tau_i\hat\gamma_X$ and $u_1\hat\gamma_X\hat\tau_i$ with $X$ of position greater than $i$, using Lemma~\ref{PreCanonicalForm}, $\hat\gamma_X\hat\tau_iu_2w_{_{>i}}$ and $\hat\tau_i\hat\gamma_Xu_2w_{_{>i}}$ are equal to a linear combination of pre-admissible words with atoms of position at least $i$. Let $v_i\dots v_r$ be one of these words. Then $w_{_{<i}}u_1v_i\dots v_r$ has less atoms of position $i$ than $w$, hence, by the induction hypothesis, it is equal to a linear combination of admissible words.
\end{itemize}
This proves that the algebra $\A$ is generated as $k$-module in degree $n$ by $A_n$.

Secondly, to prove the $k$-linear independence of the words in $A_n$, we construct $k$-algebra homomorphism $\Theta$ from $\A$ onto $\Ee$, by setting:
$\Theta(\hat\tau_i):=\tau_{\varphi_i}$, $\Theta(\hat\gamma_X):=\gamma_X$. This is indeed a surjective algebra homomorphism given that the relations in $\A$ between the $k$-generators are all satisfied by their images through $\Theta$ inside $\Ext(S_\un,S_\un)$ as showed in Lemma~\ref{TauLinRel}, Lemma~\ref{TauSqIsZero}, Lemma~\ref{GammaCommTau} and Lemma~\ref{SumGammaCommTau} in this paper, and in \cite[Proposition 12.9]{cohocplx}.
\findemo 

\npar This gives a presentation by generators and relations of $\Ee$ as stated in Conjecture~\ref{conjecture}.

\begin{mth}{Theorem}\label{PresentationExt}
The morphism $\Theta$ is an algebra isomorphism from $\A$ to $\Ext^*_{\comack_k(G)}(S_\un^G,S_\un^G)$.
\end{mth}
\pf
$\Theta:\A\to\Ext^*_{\comack_k(G)}(S_\un^G,S_\un^G)$ is surjective and, in every degree $n$, the algebra $\Ext^n_{\comack_k(G)}(S_\un^G,S_\un^G)$ has at least the dimension over $k$ of $A^n$. Hence $\Ext^n_{\comack_k(G)}(S_\un^G,S_\un^G)$ and $A^n$ have the same dimension over $k$ and $\Theta$ is an isomorphism.
\findemo

\appendix
%%%%%%%%%%%%%%%%%%%%%%%%%%%%%%%%%%%%%%%%%%%%%%%%%%%%%%%%%%%%%
\section{Arithmetics of extensions in an abelian category}%%%
%%%%%%%%%%%%%%%%%%%%%%%%%%%%%%%%%%%%%%%%%%%%%%%%%%%%%%%%%%%%%
This Appendix contain a series of classical results on computations in the graded algebra of extensions in an abelian category. We present the general framework and include the proofs of those results.

Let $\Aa$ be an abelian category. Consider following exact sequences in $\Aa$ 
$$\xymatrix{
0\ar[r]&X\ar[r]&A_1\ar[r]&A_2\ar[r]&\dots\ar[r]&A_n\ar[r]&Y\ar[r]&0}
$$
and
$$\xymatrix{
0\ar[r]&X\ar[r]&B_1\ar[r]&B_2\ar[r]&\dots\ar[r]&B_n\ar[r]&Y\ar[r]&0}
$$
representing the elements $\eta,\xi\in\Ext^n_\Aa(Y,X)$. Then $\eta+\xi$ is represented by the exact sequence
$$\xymatrix@C=3.5ex{
0\ar[r]&X\ar[r]&C_1\ar[r]&A_2\oplus B_2\ar[r]&\dots\ar[r]&A_{n-1}\oplus B_{n-1}\ar[r]&C_n\ar[r]&Y\ar[r]&0}
$$
where $C_1$, respectively $C_n$, is pushout, respectively pullback, of the following diagrams. 
$$\xymatrix{
X\oplus X\ar[r]\ar[d]_{\pi_1+\pi_2}&A_1\oplus B_1\ar[d]\\
X\ar[r]&C_1}\,,\qquad\xymatrix{
A_n\oplus B_n\ar[r]&Y\oplus Y\\
C_n\ar[r]\ar[u]&Y\ar[u]_{\Delta}}\,.
$$
When $n=1$ there is an ambiguity of this construction, given by the order in which we take the pullback and the pushout. In fact, both choices lead to equivalent exact sequences. To show this we need the following technical lemma in abelian categories. Recall that if $X$ is the pullback of the diagram $$\xymatrix{X\ar[r]\ar[d]&Y\ar[d]^f\\Z\ar[r]^g&T}$$ 
one also says that the square $(X,Y,Z,T)$ is cartesian and we have that $X$ is isomorphic to the kernel of $\xymatrix{Y\oplus Z\ar[r]^(.6){f-g}&T}$. Analogously, if $T$ is the pushout of the diagram
$$\xymatrix{X\ar[r]^s\ar[d]^t&Y\ar[d]\\Z\ar[r]&T}$$
one also says that the square $(X,Y,Z,T)$ is co-cartesian and and we have that $T$ is isomorphic to the cokernel of $\xymatrix{X\ar[r]^(.4){s\oplus t}&Y\oplus Z}$. 
%\vspace{2ex}
\pagebreak[4]
\begin{mth}{Lemma} Let $\Aa$ be an abelian category, and let
$$\xymatrix{
&0\ar[rr]&&A\ar[dl]_-a\ar[rr]^-{\alpha}&&B\ar[dl]_-b\ar[rr]^-{\beta}&&C\ar[dl]_-c\ar[rr]&&0\\
0\ar[rr]&&D\ar[rr]^-{\delta\ \ \ \ }&&E\ar[rr]^-{\varepsilon\ \ \ \ }&&F\ar[rr]&&0&\\
&0\ar'[r][rr]&&L\ar'[u][uu]^(-.4){\varphi}\ar[dl]_-l\ar'[r][rr]^(-.4){\lambda}&&M\ar@{-->}[dl]_-m\ar'[u][uu]^(-.4){\psi}\ar'[r][rr]^(-.4){\mu}&&N\ar'[u][uu]^(-.4){\chi}\ar[dl]_-n\ar[rr]&&0\\
0\ar[rr]&&P\ar[uu]^(.3)\tau\ar[rr]^-{\pi}&&Q\ar[uu]^(.3)\sigma\ar[rr]^-{\theta}&&R\ar[uu]^(.3)\rho\ar[rr]&&0&\\
}$$
be a commutative diagram with exact rows, such that the square $(A,B,D,E)$ is co-cartesian, $(Q,R,E,F)$ and $(M,N,B,C)$ are cartesian, and the morphism $n: N\to R$ is an isomorphism. \par
Then there exists a unique morphism $m:M\to L$ such that $n\mu=\theta m$ and $b\psi=\sigma m$. Moreover, the square $(L,M,P,Q)$ is commutative, and co-cartesian, and the morphisms $\tau, \varphi$ and $c$ are isomorphisms.
\end{mth}

\pf
Via the universal properties of pullback/pushout and using that the horizontal sequences are exact, one gets that $\varphi$, $\tau$ and $c$ are isomorphisms, $(A,B,D,E)$ is cartesian and, $(Q,R,E,F)$ and $(M,N,B,C)$ are co-cartesian. Now set $u=b\psi$ and $v=n\mu$. Then
$$\varepsilon u=\varepsilon b\psi=c\beta\psi=c\chi\mu=\rho n\mu=\rho v\mpoint$$
Since the square $(Q,R,E,F)$ is cartesian, there exists a unique morphism $m:M\to Q$ such that
$$u=\sigma m\;\;\hbox{and}\;\;v=\theta m\mpoint$$
Now setting $s=m\lambda$ and $t=\pi l$, we have
$$\sigma s=\sigma m\lambda=u\lambda=b\psi\lambda=b\alpha\varphi=\delta a \varphi=\delta\tau l=\sigma\pi l=\sigma t\mpoint$$
Moreover
$$\theta s=\theta m \lambda=v\lambda=n\mu\lambda=0\mvirg$$
and
$$\theta t=\theta \pi l=0\mpoint$$
Thus $\theta s=\theta t$, and $\sigma s=\sigma t$. As the square $(Q,R,E,F)$ is cartesian, this implies $s=t$, hence the square $(L,M,P,Q)$ is commutative.\par
Now let $X$ be any object of $\mathcal{A}$, and let $f:X\to M$ and $g:X\to P$ be morphisms such that $\pi g=m f$. Setting $w=\psi f$ and $r=\tau g$, we have that
$$b w=b\psi f=\sigma m f=\sigma\pi g=\delta\tau g=\delta r\mpoint$$
As $(A,B,D,E)$ is cartesian, there is a unique morphism $y:X\to A$ such that
$$w=\alpha y\;\;\hbox{and}\;\;r=a y\mpoint$$
Setting $z=\varphi^{-1}y: X\to L$, we have that
$$lz=l\varphi^{-1}y=\tau^{-1}ay=\tau^{-1}r=g\mpoint$$
Similarly
$$\psi\lambda z=\psi\lambda\varphi^{-1}y=\alpha\varphi\varphi^{-1}y=\alpha y=w=\psi f\mpoint$$

Moreover
$$\mu \lambda z=0\mvirg$$
and
$$n\mu f=\theta m f=\theta \pi g=0\mpoint$$
On the other hand
$$\chi\mu f=\beta\psi f=\beta\alpha y=0\mpoint$$
As the morphism $(\chi,n):N\to C\oplus R$ is an monomorphism, it follows that $\mu f=0=\mu\lambda z$. Since moreover $\psi\lambda z=\psi f$, and since $(M,N,B,C)$ is cartesian, it follows that $\lambda z=f$.\par
If there is another morphism $z':X\to L$ such that $l z'=g$ and $\lambda z'=f$,
then the morphism $x=z- z'$ is such that $lx=0$ and $\lambda x=0$. Then
$$\alpha\varphi x=\psi\lambda x=0\mvirg$$
and
$$a\varphi x=\tau l x=0\mpoint$$
Since $(A,B,D,E)$ is cartesian, it follows that $\varphi x=0$, hence $x=0$.\par
This shows that the square $(L,M,P,Q)$ is cartesian. Thus, to show that this square is also co-cartesian it is enough to have that $\xymatrix{P\oplus M\ar[r]^(.6){m+\pi}&Q}$ is an epimorphism. To prove this, suppose that we have $Y\in\mathcal{A}$ and $f: Q\to Y$ such that $f\pi=0$ and $fm =0$. Then $f$ factors through the cokernel of $\pi$ so there exist $g:R\to Y$ such that $f=g\theta$. Then $0=fm=g\theta m=g n \mu$.
As $\mu$ is an epimorphism and $n$ is an isomorphism, the morphism $n\mu$ is also an epimorphism, so $g=0$ and, thus $f=0$.\findemo 
When considering $\mathcal{A}$ as a module category we have an explicit description of the sum in $\Ext^1_\Aa(Y,X)$.
\begin{mth}{Lemma}\label{CompExt1}
Let $\eta,\xi\in\Ext^1_\Aa(Y,X)$.
When constructing a representative for $\eta+\xi$ both choices for the order pullback, pushout lead to the same construction.
Moreover, if $\xymatrix{0\ar[r]&X\ar[r]^i&A\ar[r]^s&Y\ar[r]&0}$ is a representative for $\eta$ and $\xymatrix{0\ar[r]&X\ar[r]^j&B\ar[r]^t&Y\ar[r]&0}$ is a representative of $\xi$. Then a representative for $\eta+\xi$ is given by $\xymatrix{0\ar[r]&X\ar[r]&C\ar[r]&Y\ar[r]&0}$ with 
$$C\simeq\{(a,b,x)\in A\oplus B\oplus X|s(a)=t(b)\}/\{(i(x_1),j(x_2),-x_1-x_2)|x_1,x_2\in X\}\,.$$
\end{mth}
\pf
To simplify notation, in the proof will be understood that $a,b,x,y$ run through $A,B,X$, respectively $Y$. Also let $$I=\{(i(x_1),j(x_2),-x_1-x_2)|x_1,x_2\in X\}\,.$$
Consider $\xymatrix{0\ar[r]&X\ar[r]&C_{ob}\ar[r]&Y\ar[r]&0}$ be the representative for $\eta+\xi$ given by taking first the pushout and then the pullback. We have
$$\xymatrix{
X\oplus X\ar[r]^{i\oplus j}\ar[d]_{\pi_1+\pi_2}&A\oplus B\ar[d]\\
X\ar[r]&C'}\,,\qquad\xymatrix{
C'\ar[r]^(0.4){\widetilde{s\oplus t}}&Y\oplus Y\\
C_{ob}\ar[r]\ar[u]&Y\ar[u]_{\Delta}}\,.
$$
with $C'\simeq A\oplus B\oplus X/I$ and, thus, 
$$C_{ob}\simeq\{(c,y)|c\in C', y\in Y, \widetilde{s\oplus t}=\Delta(y)\}\simeq\{(a,b,x)|s(a)=t(b)\}/I\,.$$
Similarly, consider $\xymatrix{0\ar[r]&X\ar[r]&C_{bo}\ar[r]&Y\ar[r]&0}$ be the representative for $\eta+\xi$ given by taking first the pullback and then pushout. We have
$$\xymatrix{
A\oplus B\ar[r]^{s\oplus t}&Y\oplus Y\\
C''\ar[r]\ar[u]&Y\ar[u]_{\Delta}}\,,\qquad\xymatrix{
X\oplus X\ar[r]^(0.6){\widetilde{i\oplus j}}\ar[d]_{\pi_1+\pi_2}&C''\ar[d]\\
X\ar[r]&C_{bo}}\,.
$$
with $C''\simeq\{(a,b,y)|s(a)=y=t(b)\}\simeq\{(a,b)|s(a)=t(b)\}$ and, thus, $C_{bo}\simeq\{(a,b,x)|s(a)=t(b)\}/I$. Hence we obtain $C_{ob}=C_{bo}$ and there is no ambiguity in the construction of a representative of $\eta+\xi$.
\findemo
We introduce here a notation we extensively use in the paper. Let $\xymatrix{X\ar[r]^f&A\ar[r]^g&Y}$ be an exact sequence at $A$. Then we write $A\simeq\bn{\Coker\,f}{\Ker\,g}$. In particular, we have a short exact sequence
$$\xymatrix{
0\ar[r]&X\ar[r]&\bn{Y}{X}\ar[r]&Y\ar[r]&0}\,.
$$
and a representative in $\Ext^2_\Aa(Y,X)$ can be written as
$$\xymatrix{
0\ar[r]&X\ar[r]&\bn{Z}{X}\ar[r]&\bn{Y}{Z}\ar[r]&Y\ar[r]&0}\,.
$$
With this notation we have an easy way to compute the exact sequence representing the sum in $\Ext^2_\Aa(Y,X)$. The computation is an easy consequence of the following technical lemma whose proofs are left to the reader.
\begin{mth}{Lemma}\label{ExactCoCartesian} Let $\mathcal{A}$ be an abelian category, and let
$$\xymatrix{
0\ar[r]&X\oplus X\ar[d]_-a\ar[r]^-{\alpha}&{\bn{A}{X}\oplus\bn{B}{X}}\ar[d]_-b\ar[r]^-{\beta}&A\oplus B\ar[d]_-c\ar[r]&0\\
0\ar[r]&X\ar[r]^-{\delta}&C\ar[r]^-{\varepsilon}&A\oplus B\ar[r]&0\\
}$$
be a commutative diagram with exact rows, such that the square on the left is co-cartesian. Then $c$ is an isomorphism and, hence, $C\simeq\bn{A\oplus B}{X}$. 
\end{mth}
and its dual
\begin{mth}{Lemma}\label{ExactCartesian} Let $\mathcal{A}$ be an abelian category, and let
$$\xymatrix{
0\ar[r]&A\oplus B\ar[r]^-{\alpha}&{\bn{Y}{A}\oplus\bn{Y}{B}}\ar[r]^-{\beta}&Y\oplus Y\ar[r]&0\\
0\ar[r]&A\oplus B\ar[r]^-{\delta}\ar[u]_-f&C'\ar[r]^-{\varepsilon}\ar[u]_-g&Y\ar[u]_-h\ar[r]&0\\
}$$
be a commutative diagram with exact rows, such that the square on the right is cartesian. Then $f$ is an isomorphism and, hence, $C'\simeq\bn{Y}{A\oplus B}$. 
\end{mth}
Remark that, in the above notation, $\bn{A\oplus B}{X}$ has $\bn{A}{X}$ and $\bn{B}{X}$ as subobjects and $\bn{Y}{A\oplus C}$ has $\bn{Y}{A}$ and $\bn{Y}{C}$ as factors. We conclude with the computation of the exact sequence representing the sum in $\Ext^2_\Aa(Y,X)$.
\begin{mth}{Lemma}\label{CompExt2}
Let $\eta,\xi\in\Ext^2_\Aa(Y,X)$ be represented by
$$\xymatrix{
0\ar[r]&X\ar[r]&\bn{Z_1}{X}\ar[r]&\bn{Y}{Z_1}\ar[r]&Y\ar[r]&0}\,,
$$
respectively by
$$\xymatrix{
0\ar[r]&X\ar[r]&\bn{Z_2}{X}\ar[r]&\bn{Y}{Z_2}\ar[r]&Y\ar[r]&0}\,.
$$
Then $\eta+\xi$ is represented by
$$\xymatrix{
0\ar[r]&X\ar[r]&\bn{Z_1\oplus Z_2}{X}\ar[r]&\bn{Y}{Z_1\oplus Z_2}\ar[r]&Y\ar[r]&0}\,.
$$
\end{mth}
\pf
Apply Lemma~\ref{ExactCoCartesian} and Lemma~\ref{ExactCartesian}.
\findemo
%\bibliographystyle{abbrv}
%\bibliography{mabib}

Serge Bouc - CNRS-LAMFA - Universit\'e de Picardie - 33, rue St Leu - 80039 - Amiens - France. \\
Email~: {\tt serge.bouc@u-picardie.fr}\spn
Radu Stancu - LAMFA - Universit\'e de Picardie - 33, rue St Leu - 80039 - Amiens - France. \\
Email~: {\tt radu.stancu@u-picardie.fr}

\end{document}